\documentclass{amsart}
\usepackage{amsmath, amssymb, amsthm,amsfonts}
\usepackage{enumitem}
\usepackage{graphicx,color}
\usepackage{graphics}
\usepackage{comment}
\usepackage[a4paper, top = 1.2in, bottom = 1.2in, left = 1.2in, right = 1.2in]{geometry}
\usepackage[OT2,T1]{fontenc}
\DeclareSymbolFont{cyrletters}{OT2}{wncyr}{m}{n}
\DeclareMathSymbol{\Sha}{\mathalpha}{cyrletters}{"58}

\DeclareMathOperator{\LCM}{lcm}

\DeclareMathOperator{\SL}{SL}
\DeclareMathOperator{\GL}{GL}

\DeclareMathOperator{\cond}{Cond}

\newcommand{\Nsc}{\N^{\mathrm{sc}}}
\newcommand{\Hup}{\mathbb{H}}

\newcommand{\bchi}{\boldsymbol{\chi}}

\newcommand{\ord}{{\rm ord}}
\newcommand{\SHA}[1]{{\cyr X}(#1)}

\newcommand{\As}{\mathbb{A}_\Q^{\times}}

\newcommand{\Qs}{\mathbb{Q}^{\times}}
\newcommand{\Zs}{\mathbb{Z}^{\times}}
\newcommand{\Q}{{\mathbb Q}}
\newcommand{\Z}{{\mathbb Z}}
\newcommand{\N}{{\mathbb N}}
\newcommand{\C}{{\mathbb C}}
\newcommand{\R}{{\mathbb R}}
\newcommand{\F}{{\mathbb F}}

\newcommand{\OO}{{\mathcal O}}

\newcommand{\kro}[2]{\left( \frac{#1}{#2} \right) }

            \DeclareFontFamily{U}{wncy}{} 
            \DeclareFontShape{U}{wncy}{m}{n}{%
               <5>wncyr5%
               <6>wncyr6%
               <7>wncyr7%
               <8>wncyr8%
               <9>wncyr9%
               <10>wncyr10%
               <11>wncyr10%
               <12>wncyr6%
               <14>wncyr7%
               <17>wncyr8%
               <20>wncyr10%
               <25>wncyr10}{} 
\DeclareMathAlphabet{\cyr}{U}{wncy}{m}{n}

\begin {document}

\newtheorem{thm}{Theorem}

\newtheorem{lem}{Lemma}[section]
\newtheorem{prop}[lem]{Proposition}

\newtheorem{cor}[lem]{Corollary}

\theoremstyle{definition}

\newtheorem{ex}{Example}

\theoremstyle{remark}

\newtheorem*{ack}{Acknowledgement}

\title[Explicit application of Waldspurger's Theorem]{
Explicit application of Waldspurger's Theorem
}

\author{Soma Purkait}
\address{Mathematics Institute\\
	University of Warwick\\
	Coventry\\
	CV4 7AL \\
	United Kingdom}

\email{S.Purkait@warwick.ac.uk}
\date{\today}
\thanks{The author is supported by a Warwick Postgraduate Research Scholarship
and an EPSRC Delivering Impact Award}


\begin{abstract}
For a given cusp form $\phi$ of even integral weight satisfying certain hypotheses, Waldspurger's Theorem
relates the critical value of the $\mathrm{L}$-function of the $n$-th quadratic twist of $\phi$ to
the $n$-th coefficient
of a certain modular form of half-integral weight. Waldspurger's recipes for these modular forms of half-integral
weight are far from being explicit. 
In particular, they are expressed in the language of automorphic representations and Hecke characters.
We translate these recipes into congruence conditions involving 
easily computable values of Dirichlet characters. We illustrate the practicality of our \lq simplified
Waldspurger\rq\ by giving several examples.
\end{abstract}
\maketitle

\section{Introduction}

In 1983 Tunnell \cite{Tunnell} gave a remarkable solution to the congruent number
problem, assuming the celebrated Birch and Swinnerton-Dyer Conjecture. This
ancient Diophantine question asks for the classification of congruent numbers,
those positive integers which are the areas of right-angled triangles whose sides
are rational numbers. 
For positive $n$, write $E_n : y^2=x^3-n^2 x$;  note
that $E_n$ is the $n$-th quadratic twist of $E_1$. It is straightforward
to show that $n$ is a congruent number if and only if the elliptic curve
$E_n/\Q$ has positive rank. Tunnell expresses the critical value
of the $\mathrm{L}$-function of $E_n$ in terms of coefficients of certain modular forms of weight $3/2$.
These modular forms are in turn written in terms of theta series of ternary
quadratic forms. Applying the conjecture of Birch and Swinnerton-Dyer, Tunnell is then able
to give a simple and elegant criterion for $n$ to be a congruent number.

Tunnell's Theorem is a highly non-trivial consequence of a theorem of Waldspurger \cite{Waldspurger}.
For a given cusp form $\phi$ of even integral weight satisfying certain hypotheses, Waldspurger's Theorem
relates the critical value of the $\mathrm{L}$-function of the $n$-th quadratic twist of $\phi$ to
the $n$-th coefficient
of a certain modular form of half-integral weight. Waldspurger's recipes for these modular forms of half-integral
weight are far from being explicit. 
In particular, they are expressed in the language of automorphic representations and Hecke characters.
We translate these recipes into congruence conditions involving 
easily computable values of Dirichlet characters. We illustrate the practicality of our \lq simplified
Waldspurger\rq\ by giving several Tunnell-like examples, of which the following is the simplest.

\begin{prop}\label{prop:50}
Let $E$ be the elliptic curve of conductor $50$ given by
\begin{equation}\label{eqn:50}
 E: Y^2 + XY +Y = X^3 +X^2 -3X +1.
\end{equation} 
Let $Q_1$, $Q_2$, $Q_3$, $Q_4$ be the following positive-definite ternary quadratic forms,
\begin{equation*}
\begin{split}
Q_1 = 25x^2+ 25 y^2+z^2,  \qquad Q_2 = 14x^2+9y^2+6z^2+4yz+6xz+2xy,\\ 
Q_3 = 25 x^2+13y^2+2z^2+2yz, \qquad Q_4 = 17x^2+17y^2+3z^2-2yz-2xz+16xy. 
\end{split}
\end{equation*}
Let $n$ be a positive square-free number such that $5 \nmid n$. Then,
\[
\mathrm{L}(E_{-n},1) = \frac{\mathrm{L}(E_{-1},1)}{\sqrt{n}}
\cdot c_n^2 \]
where $E_{-n}$ is the $-n$-th quadratic twist of $E$ and 
\[ c_n = \sum_{i=1}^4 \frac{(-1)^{i-1}}{2} \cdot \# \{ (x,y,z) :Q_i(x,y,z) = n\}. \] 
\end{prop}
For the elliptic curve $E$ in \eqref{eqn:50} and $n$ a positive square-free integer such
that $5 \nmid n$, we give a similar formula for $\mathrm{L}(E_n,1)$ that involves
$38$ quadratic forms.

We would like to note here that all the examples we consider deal with newforms whose levels are 
neither odd nor square-free. In fact for newforms $\phi$ of weight $k-1$ and odd and square-free level $N$
with $\mathrm{L}(\phi,1) \ne 0$ and $k \equiv 3 \pmod{4}$ there is an explicit recipe by Bocherer 
and Schulze-Pillot ~\cite{B-S2} for constructing a modular form of weight $k/2$, level $4N$ that is Shimura 
equivalent to $\phi$. Their method uses generalized theta series and the Eichler correspondence 
with automorphic forms on quaternion algebras. In particular they show that given a rational 
elliptic curve $E$ of odd and square-free conductor, an inverse Shimura lift of $\phi_E$ 
(the newform corresponding to $E$) comes from ternary quadratic forms if and only if 
$\mathrm{L}(E, 1) \ne 0$. We note that the form they construct belongs to the Kohnen subspace~\cite{Kohnen}. It 
follows by Waldspurger~\cite{Waldspurger} that in these cases the space of Shimura 
equivalent forms at level $4N$ is two-dimensional. In a recent paper Hamieh~\cite{Hamieh} used~\cite{B-S2} and 
Waldspurger's recipe to compute a basis for this two-dimensional space.

We would also like to mention the work of Shin-ichi Yoshida~\cite{Yo} in which he considers 
$2\pi/3$ and $\pi/3$-congruent number problems and uses Waldspurger's result (Corollary~\ref{cor:wald} below) 
to give a Tunnell-like criterion for a square-free 
number in certain congruence classes to be $2\pi/3$ and $\pi/3$-congruent. 

The paper is arranged as follows. In Section~\ref{section:Shimura} we review 
Shimura's decomposition of the space of cusp forms
of a certain level and half-integral weight into certain subspaces
appearing in Waldspurger's Theorem.
In Section~\ref{section:Heckechar} we review the 
correspondence between Dirichlet characters
and Hecke characters and we prove a result that allows us to
evaluate components of a Hecke character corresponding to a given Dirichlet character. 
Next, in Section~\ref{sec:aut} we 
review the correspondence between modular forms
of even integral weight and automorphic representations
and prove a result needed for simplifying the hypotheses
of Waldspurger's Theorem. 
In Section~\ref{section:waldspurger'sthm} we state Waldspurger's Theorem
in simplified form.
To apply Waldspurger's Theorem in conjunction with the
Birch and Swinnerton-Dyer Conjectures it is convenient
to express the period of the $n$-th twist of a given
elliptic curve in terms of the period of the elliptic curve itself.
We do this in Section~\ref{sec:period}. 
To apply Waldspurger's recipes we need to be able to
answer questions of the following form: for a given
cusp form of half-integral weight $f=\sum a_n q^n$,
and positive integers $a$, $M$, is $a_n=0$ for all
$n \equiv a \pmod{M}$? We give an algorithm for answering
this question in Section~\ref{section:coeffmodulon}.
Finally, Section~\ref{section:appl}  is devoted to
extensive examples which combine our algorithm~\cite{SomaI} for computing 
Shimura decomposition with Waldspurger's Theorem as made explicit in this paper.

\begin{ack} 
This article is a part of my thesis and I would like to thank my advisor 
Professor Samir Siksek for suggesting this problem and for several helpful 
discussions. I am grateful to Professor John Cremona and Dr. Neil Dummigan for 
careful reading of my thesis and for providing valuable insights and suggestions.
I am grateful to the referee for providing many helpful remarks and suggestions and 
pointing out many useful references.
\end{ack} 

\section{Shimura Decomposition}\label{section:Shimura}
Let $k\geq 3$ be an odd integer and $N$ a positive integer such that $4 \mid N$.
Let $\chi$ be an even Dirichlet character modulo $N$.
We denote by $S_{k/2}(N,\chi)$ the space of cusp forms of weight $k/2$, level $N$ and 
character $\chi$. Let $S_{k/2}^0(N,\chi)$ be the subspace of $S_{k/2}(N,\chi)$ spanned by 
single-variable\footnote{The term ``single-variable theta series'' refers to the theta series 
of weights $1/2$ and $3/2$ that come from a quadratic form of one variable and are of the form 
$\sum_{n= -\infty}^{\infty} \psi(n)n^\nu q^{n^2}$, 
where $\nu \in \{0,1\}$ and $\psi$ is a Dirichlet character such that $\psi(-1) = (-1)^\nu$.} theta series when $k=3$; 
for $k\geq 5$, we define $S_{k/2}^0(N,\chi)=0$.
More precisely, a generating set for $S_{3/2}^0(N,\chi)$ is given by
\begin{equation*}
\begin{split}
S=\{\ \sum_{m=1}^{\infty} \psi(m)mq^{tm^2} :\ 4 r_{\psi}^2 t \mid N\ \text{and}\ & \text{$\psi$ is 
a primitive odd character of} \\
& \text{conductor $r_\psi$ such that}\ \chi = \kro{-4t}{.}\psi \ \},
\end{split}
\end{equation*}
which in fact constitutes a basis for $S_{3/2}^0(N,\chi)$, as shown in ~\cite{SomaI}. 
The interesting part (from the point-of-view of 
Waldspurger's Theorem) of the space $S_{k/2}(N,\chi)$ is the orthogonal complement of $S_{k/2}^0(N,\chi)$
with respect to the Petersson inner product, denoted by $S_{k/2}^\prime(N,\chi)$. 

In his thesis Basmaji~\cite{Basmaji} gave an algorithm for computing a basis for the space of 
half-integral weight modular forms of level divisible by $16$. The main idea of the algorithm is to use theta series 
$\Theta = \sum_{n= -\infty}^{\infty} q^{n^2}$, 
$\Theta_1 =  \frac{\Theta - V(4)\Theta}{2}$ and the following embedding,
\[
\varphi : S_{k/2}(N, \chi) \rightarrow {S \times S},
\qquad \qquad
f \mapsto (f\Theta,f\Theta_1),
\]
where $S= S_{\frac{k+1}{2}}\left(N,\ \chi\cdot\chi_{-1}^\frac{k+1}{2}\right)$ and $V$ is the usual shift operator. This idea 
has been generalized by Steve Donnelly for level divisible by $4$ and is implemented in {\tt MAGMA}.

Let $N^\prime=N/2$. For $M \mid N^\prime$ such that $\cond(\chi^2) \mid M$
and a newform $\phi \in S_{k-1}^{\mathrm{new}}(M,\chi^2)$ Shimura defines
\[
S_{k/2}(N,\chi,\phi)=
\{
f \in S_{k/2}^{\prime}(N,\chi) : 
\text{$T_{p^2}(f)=\lambda_p(\phi) f$ for almost all $p \nmid N$}\},
\]
here $T_{p} (\phi) =\lambda_p(\phi) \phi$, and gives the following decomposition theorem \cite{ShimuraII}: 
\begin{thm}(Shimura)
$S_{k/2}^{\prime}(N,\chi)=\bigoplus_\phi S_{k/2}(N,\chi,\phi)$ where $\phi$
runs through all newforms $\phi \in S_{k-1}^{\mathrm{new}}(M,\chi^2)$
with $M \mid N^\prime$ and $\cond(\chi^2) \mid M$.
\end{thm}

We point out that the summands $S_{k/2}(N,\chi,\phi)$ occur in Waldspurger's Theorem
and their computation is necessary for explicit applications of that theorem. 
However the above theorem is not suitable for computation since for any particular prime 
$p \nmid N$, we do not know if it is included or excluded in the \lq almost all\rq\ condition. 
In ~\cite{SomaI} we proved the above theorem with a more precise definition for the spaces 
$S_{k/2}(N,\chi,\phi)$: 
\[
S_{k/2}(N,\chi,\phi)=
\{
f \in S_{k/2}^{\prime}(N,\chi) : 
\text{$T_{p^2}(f)=\lambda_p(\phi) f$ for  all $p \nmid N$}\},
\] 
whilst showing that our definition is equivalent to Shimura's definition.
We also proved the following theorem that gives an algorithm for computing the Shimura decomposition.
\begin{thm}\cite{SomaI}\label{thm:algo}
Let $\phi$ be a newform of weight $k-1$, level $M$ dividing $N^\prime$, and character $\chi^2$.
Let $p_1,\dots,p_n$ be primes not dividing $N$ satisfying the following: for every newform 
$\phi^\prime \ne \phi$ of weight $k-1$, level dividing $N^\prime$ and character $\chi^2$, there
is some $p_i$ such that $\lambda_{p_i}(\phi^\prime) \ne \lambda_{p_i}(\phi)$, where
$T_{p_i}(\phi)=\lambda_{p_i}(\phi) \cdot \phi$. Then
\[
S_{k/2}(N,\chi,\phi)=
\left\{
f \in S_{k/2}(N,\chi) \; : \;
T_{p_i^2} (f) = \lambda_{p_i} (\phi) f  \quad \text{for $i=1,\dots,n$} 
\right\}.
\]
\end{thm}

\section{Correspondence between Dirichlet Characters and Hecke
Characters on $\As/\Qs$ of Finite Order}\label{section:Heckechar}

We shall need the correspondence between Dirichlet characters
and Hecke characters on $\As/\Qs$ of finite order. This
material is in Tate's thesis \cite{CF}, but 
we found the presentation in \cite[Section 3.1]{Bump} more useful.

\begin{prop} Let $\bchi=(\chi_p)$
be a character on $\As$. Then
there exists a finite set $S$ of 
places, including all the 
Archimedean ones, such that if 
$p \notin S$, then $\chi_p$ is
trivial on the unit group $\Z_p^{\times}$.
\end{prop}

Recall that if $\chi_p$ is trivial on the unit group $\Z_p^{\times}$, then $\chi_p$ is unramified.
Thus by the above proposition, $\chi_p$ is unramified for all but finitely many $p$.

\begin{thm} (\cite[Proposition 3.1.2]{Bump}) \label{thm:correspond}
Suppose $\boldsymbol{\chi}=(\chi_p)$ is a character of 
finite order on $\As/\Qs$.
There exists an integer $N$ whose
prime divisors are precisely the
non-Archimedean primes $p$ such that
$\chi_p$ is ramified, and a
primitive Dirichlet character $\chi$
modulo $N$ such that
if $p\nmid N$ is non-Archimedean then
$\chi(p)= \chi_p(p)$.
This correspondence $\bchi \mapsto \chi$
is a bijection between characters
of finite order of $\As/\Qs$
and the primitive Dirichlet characters.
\end{thm}

In our work, we shall need to start
with a Dirichlet character $\chi$ of modulus $N$
and then do computations with
the corresponding Hecke 
character $\bchi$. We collect
here some facts that will help
us with these computations.

\begin{lem}\label{lem:charprop}
We keep the notation of Theorem~\ref{thm:correspond}.
\begin{enumerate}
\item[(i)] For any $\alpha \in \Qs$, $\prod \chi_p(\alpha)=1$ where the product is taken over all places.
\item[(ii)] Suppose $p=\infty$ and $\alpha \in \Q_\infty^{\times}=\R^\times$.
Then $\chi_\infty(\alpha)=1$ if $\alpha>0$, or if $\chi$ has odd order.
\item[(iii)] Let $p$ be a non-Archimedean prime such that $p \mid N$ and $\alpha$, $\beta \in \Z_p$ be non-zero.
Suppose that $\beta \equiv \alpha \pmod{\alpha N \Z_p}$. Then
$\chi_p(\beta)=\chi_p(\alpha)$.
\item[(iv)] Let $p$ be non-Archimedean such that $p \nmid N$. Then $\chi_p$ is unramified.
\end{enumerate}
\end{lem}
\begin{prop} \label{prop:ev}
Let $\chi$ be a Dirichlet character modulo $N$ (not necessarily primitive)
and let $\bchi=(\chi_p)$ be the corresponding character on $\As/\Qs$.
Let $a \in \Z$ be non-zero. For a prime $q$, let $\nu_q(a)$ denotes the 
highest power of $q$ that divides $a$.
\begin{enumerate}
\item[(a)] If $q \nmid N$ then $\chi_q(a)=\chi(q)^r$ where $r=\nu_q(a)$.
\item[(b)] Suppose $q$ divides $N$ and let $q_1,\dotsc,q_r$ be the other
primes dividing $N$. Let $b$ be a {\bf positive} integer satisfying
\[
b \equiv \begin{cases} a \pmod {aN\Z_q} \\
1 \pmod {N\Z_{q_i}} & i =1,\dots,r;
\end{cases}
\]
such  $b$ can easily be constructed by the Chinese Remainder Theorem.
Write
\[
b=q^{\nu_q(a)} \prod_{j=1}^s \ell_j^{\beta_j}
\]
where the $\ell_j$ are distinct primes. 
Then
\[
\chi_q(a)=\prod_{j=1}^s \chi(\ell_j)^{-\beta_j}.
\]
\end{enumerate}
\end{prop}
\begin{proof}
Let $N^\prime$ be the conductor of $\chi$ and note that $N^\prime \mid N$.
Now if $q \nmid N$ then, $\chi_q$ is unramified. Write
$a=q^r a^\prime$ where $q \nmid a^\prime$. Then $a^\prime \in \Z_q^{\times}$. Thus by definition
of unramified, $\chi_q(a^\prime)=1$. Moreover, from Theorem~\ref{thm:correspond},
$\chi_q(q)=\chi(q)$. This proves (a).

Now suppose $q \mid N$ and let $q_1,\dotsc,q_r$ be the other
primes dividing $N$. Let $b$ be as in the proposition. Since $N^\prime \mid N$,
we have
\[
b \equiv \begin{cases} a \pmod { aN^\prime \Z_q} \\
1 \pmod {N^\prime \Z_{q_i}} & i =1,\dots,r.
\end{cases}
\]
By Lemma~\ref{lem:charprop},
$\chi_q(b)=\chi_q(a)$, and $\chi_{q_i}(b)=1$ for $i=1,\dots,r$.
Now 
\begin{equation*}
\begin{split}
\chi_q(a) &= \chi_q(b)\\
&= \prod_{p \neq q} \chi_p(b)^{-1} \qquad \text{by (i) of Lemma~\ref{lem:charprop},} \\
&= \prod_{p \nmid N} \chi_p(b)^{-1} \qquad \text{since $\chi_{q_i}(b)=1$,}\\
&=\prod_{j=1}^s \chi(\ell_j)^{-\beta_j} \qquad \text{using part (a)}.
\end{split}
\end{equation*}
This completes the proof.
\end{proof}
\section{Local components of the automorphic representations associated to modular forms of even integer weight} \label{sec:aut}

Let $k$ be a positive odd integer with $k\geq 3$. Let $\phi = \sum_{n=1}^{\infty} a_nq^n 
\in S_{k-1}^{\mathrm{new}}(N,\chi)$ be 
a newform of weight $k-1$, level $N$ and character $\chi$.

We can associate to $\phi$ an automorphic representation $\rho$. Let $\rho_p$ be the
local component of $\rho$ at a prime $p$.

If $\phi=\sum_{n=1}^\infty a_n q^n$ is an eigenform,
then we define its twist by a character $\mu$ to
be the modular form $\phi_\mu=\sum_{n=1}^\infty a_n \mu(n) q^n$.

Waldspurger works with the following different definition of twist:
Let $\phi$ be a newform of weight $k-1$ and character
$\chi$. Let $\mu$
a Dirichlet character. We denote by $\phi \otimes \mu$ the (unique) newform
of weight $k-1$ with character  $\chi \mu^2$ satisfying 
$\lambda_p(\phi \otimes \mu)=\mu(p) \lambda_p(\phi)$ for almost all primes $p$,
where $\lambda_p$ is the eigenvalue under $T_p$.

Now fix a prime number $p$. Let $\xi_p$ be the set of primitive Dirichlet characters with 
$p$-power conductor. The following holds (see \cite[Section III]{Waldspurger}):

\begin{enumerate}
\item[(i)] $\rho_p$ is supercuspidal if and only if for all $\mu \in \xi_p$, the level of $\phi \otimes \mu$ 
is divisible by $p$ and $\lambda_p(\phi \otimes \mu)=0$.
\item[(ii)] $\rho_p$ is an irreducible principal series if and only if either 
\begin{enumerate}
\item[(a)] there exists a character $\mu$ in $\xi_p$ such that $p$ does not 
divide the level of $\phi \otimes \mu$; or
\item[(b)] there exist two distinct characters $\mu_1$, $\mu_2$ in $\xi_p$ such that $\lambda_p(\phi \otimes \mu_1) \ne 0$, 
$\lambda_p(\phi \otimes \mu_2) \ne 0$.
\end{enumerate}
\item[(iii)] $\rho_p$ is a special representation if and only if the following conditions hold: 
\begin{enumerate}
\item[(a)] for all $\mu \in \xi_p$, the level of $\phi \otimes \mu$ is divisible by $p$; 
\item[(b)] there exists a unique $\mu$ in $\xi_p$ such that $\lambda_p(\phi \otimes \mu) \ne 0$.
\end{enumerate}
\end{enumerate}

We shall need the following theorem which is extracted from the
paper of Atkin and Li \cite{Atkin-Li}.
\begin{thm} (Atkin and Li) \label{thm:AtkinLi}
Let $\phi = \sum_{n=1}^\infty a_nq^n$ be a newform of weight $k-1$, character $\chi$
and level $N$. Let $\mu$ be a primitive character of conductor
$m$. Then
\begin{enumerate}
\item[(a)] If $\gcd(m,N)=1$ then $\phi\otimes \mu=\phi_\mu$,
and it is a newform of weight $k-1$, character $\chi \mu^2$
and level $N m^2$ \cite[Introduction]{Atkin-Li}.

\item[(b)]Suppose $\mu$ is of $q$-power conductor where $q \mid N$
and write $N=q^s M$ where $q \nmid M$.
Then $\phi\otimes \mu$ is a newform of weight $k-1$, character
$\chi \mu^2$ and level $q^{s^\prime} M$ for some $s^\prime \geq 0$.
Moreover, $\lambda_p(\phi \otimes \mu)=\mu(p)\lambda_p(\phi)$
for all primes $p \nmid N$ \cite[Theorem 3.2]{Atkin-Li}.
In particular if $s=1$ and $\chi$ is trivial, then for $\mu$ with 
conductor $q^r$, $r\geq 1$, it turns out that $\phi\otimes \mu = \phi_\mu$ 
is a newform of level $q^{2r}M$ and character $\mu^2$ \cite[Corollary 4.1]{Atkin-Li}.

\item[(c)]Let $q \mid N$. Suppose $\phi$ is $q$-primitive and $a_q = 0$.
Then for all characters $\mu$ of $q$-power conductor, $\phi\otimes \mu=
\phi_\mu$ is a newform of level divisible by $N$ (Recall that $\phi$ is $q$-primitive 
if $\phi$ is not a twist of any newform of level lower than $N$ by a character of
conductor equal to some power of $q$) \cite[Proposition 4.1]{Atkin-Li}.

\item[(d)] Let $N=q^s M$ where $q \nmid M$; let $Q=q^s$. 
Let $\chi_Q$ be the $Q$-part~\footnote{
Let $\chi$ be a Dirichlet character 
with modulus $p_1^{r_1} \cdots p_n^{r^n}$ where the $p_i$
are distinct primes. Then $\chi$ can be written
uniquely as a product $\prod \chi_{p_i^{r_i}}$
where $\chi_{p_i^{r_i}}$ has modulus $p_i^{r_i}$. See \cite{Atkin-Li}.
}
of the character $\chi$.
If $s$ is odd and $\mathrm{cond}\ \chi_Q \leq \sqrt{Q}$ then $\phi$ is $q$-primitive.

Now suppose $q=2$. Then,
if $s=2$ then $\phi$ is always $2$-primitive; 
if $s$ is odd then $\phi$ is $2$-primitive if and only if $\mathrm{cond}\ \chi_Q < \sqrt{Q}$; 
if $s$ is even and $s\geq 4$ then $\phi$ is $2$-primitive if and only if $\mathrm{cond}\ \chi_Q = \sqrt{Q}$ 
\cite[Theorem 4.4]{Atkin-Li}.
\end{enumerate}
\end{thm}

We deduce the following corollaries which we will be using later.
\begin{cor}\label{cor:supercuspidal}
Let $\phi = \sum_{n=1}^{\infty} a_nq^n \in S_{k-1}^{\mathrm{new}}(N)$ be a newform with trivial character. 
Let $\rho_2$ be the local component at $2$ of the corresponding automorphic representation. Suppose either
$N$ is odd or $\nu_2(N)=1$. Then $\rho_2$ is not supercuspidal.
Further if $\nu_2(N) \geq 2$ and $\phi$ is $2$-primitive then $\rho_2$ is supercuspidal, hence   
if either $\nu_2(N)=2$ or $\nu_2(N)>1$ is odd then $\rho_2$ is supercuspidal.
\end{cor}

\begin{proof}
If $N$ is odd, take $\mu$ to be the identity character. Thus $\mu \in \xi_2$ and the level of $\phi\otimes \mu$ 
is odd and hence $\rho_2$ is not supercuspidal. If $N=2M$ such that $M$ is odd then $a_2 \ne 0$, so taking $\mu$ 
as identity character we get that $\lambda_2(\phi\otimes \mu) = a_2 \ne 0$ and thus $\rho_2$ is not supercuspidal. 

Let $\nu_2(N) \geq 2$. Then $a_2 = 0$. If $\phi$ is $2$-primitive then 
it follows using part $(c)$ of Theorem \ref{thm:AtkinLi} that for any $\mu \in \xi_2$, 
$\phi\otimes \mu = \phi_\mu$ is newform of level divisible by $2$. 
Write $T_2(\phi_\mu)= \sum_{n=1}^{\infty} b_n q^n$.
Then, 
$b_n = a_{2n}\mu(2n) + {\mu}^2(2) 2^{k-2}a_{n/2}\mu(n/2)$
for all $n$. Thus $T_2(\phi_\mu)=0$. Therefore,
 $\lambda_2(\phi\otimes \mu) = \lambda_2(\phi_\mu)= 0$ and $\rho_2$ is supercuspidal. 
The final statement is a direct application of part $(d)$ of Theorem \ref{thm:AtkinLi}. 
\end{proof}

\begin{cor}\label{cor:special}
Let $\phi$ be as in the above corollary. 
\begin{enumerate}
\item[(i)] If $N = pM$ with $M$ coprime to $p$ and $a_p \ne 0$, then $\rho_p$ is a special representation.
\item[(ii)]
If $p \nmid N$, then $\rho_p$ is an irreducible principal series.
\end{enumerate}
\end{cor}
\begin{proof}
We first prove (i). 
By part (b) of the Theorem \ref{thm:AtkinLi}, for any $\mu \in \xi_p$, the level of $\phi\otimes \mu$ is divisible by $p$. 
Further if $\mu$ is the identity character then $\lambda_p(\phi\otimes \mu) = a_p \ne 0$; we claim that this is unique such character in $\xi_p$. 
Let $\mu \in \xi_p$ be such that $\mu$ is a character of conductor $p^r$, $r\geq 1$. Then 
$\phi\otimes \mu = \phi_\mu$ is a newform in $S_{k-1}(p^{2r} M, {\mu}^2)$ such that $\lambda_p(\phi_\mu)= a_p\mu(p) = 0$ 
and hence $\lambda_p(\phi\otimes \mu)=0$. 

The proof of (ii) is obvious and does not require the condition that newform $\phi$ has trivial character.
\end{proof}

\section{Waldspurger's Theorem and Notation}\label{section:waldspurger'sthm}
In this section we will present Waldspurger's Theorem. We will introduce and simplify
the notation used in the theorem. This is needed in the following section where 
we will discuss how to use the theorem for elliptic curves and compute critical values of $\mathrm{L}$-functions
in terms of coefficients of corresponding half-integral weight forms.
An important application is the computation of
orders of the Tate-Shafarevich groups assuming the Birch and Swinnerton-Dyer Conjecture.

Let $k$ be positive integers with $k\geq 3$ odd. Let $\chi$ be an even Dirichlet character with modulus 
divisible by $4$. Fix a newform $\phi$ of level $M_\phi$ in $S_{k-1}^{\mathrm{new}}(M_{\phi},\chi^2)$. 
Let $p$ be a prime number. Let $\nu_p$ be the $p$-adic valuation on $\Q$ and $\Qs_p$.
Let $m_p = \nu_p(M_{\phi})$ and $\lambda_p$ be the Hecke eigenvalue of $\phi$ corresponding
to the Hecke operator $T_p$.  

Let $\rho$ be the automorphic representation associated to $\phi$ and $\rho_p$ be the
local component of $\rho$ at $p$. Let $S$ be the (finite) set of primes $p$ such that 
$\rho_p$ is not irreducible principal series. 
If $p \notin S$, $\rho_p$ is equivalent to $\pi(\mu_{1,p}, \mu_{2,p})$ 
where $\mu_{1,p}$ and $\mu_{2,p}$ are two continuous characters on $\Q_p$ such that 
$\mu_{1,p} \mu_{2,p} \ne \lvert \cdot \rvert^{\pm 1}$.  
Let (H1) be the following hypothesis:
\[
\mathrm{(H1)} \qquad  
\text{For all $p \notin S$, $\mu_{1,p}(-1) = \mu_{2,p}(-1) = 1$.}
\]

\begin{thm} (Flicker) \label{thm:Flicker}
There exists $N$ such that $S_{k/2}(N,\chi,\phi) \ne \{0\}$ if and only if the 
hypothesis $\mathrm{(H1)}$ holds.
\end{thm}

\begin{thm} \label{thm:Vig} (Vigneras)
Flicker's condition $\mathrm{(H1)}$ 
always holds whenever $\phi$ is a newform of even weight with trivial character.
\end{thm}
\begin{proof}
 For the proof refer to \cite{Vigneras}.
\end{proof}

From the theorems of Flicker and Vigneras we have the following easy corollary.
\begin{cor}
Let $\phi$ be a newform of weight $k-1$, level $M_\phi$ and trivial character
$\chi_{\mathrm{triv}}$. Let $\chi$ be a Dirichlet character 
satisfying  $\chi^2 = \chi_{\mathrm{triv}}$.
Then there exists some $N$ such that $S_{k/2}(N,\chi,\phi) \ne \{0\}$.
\end{cor}

Henceforth, we will always assume that $\phi$ has trivial character and $\chi$ is
quadratic, thus the conclusion of the corollary holds.
We will now introduce several pieces of notation used by Waldspurger \cite[Section VIII]{Waldspurger} 
before stating his main theorem.

Let $\chi_0$ be the Dirichlet character associated to $\chi$ given by 
\[
\chi_0(n) := \chi(n) \kro{-1}{n}^{(k-1)/2}.
\]
Let $\chi_{0,p}$ be the local component of $\chi_0$ at a prime $p$. 
For each prime $p$ we will later define a non-negative integer $\widetilde{n_p}$ 
that depends only on the local components $\rho_p$ and $\chi_{0,p}$. 
Let $\widetilde{N_{\phi}}$ be given by 
\[
\widetilde{N_{\phi}} := \prod_p{p^{\widetilde{n_p}}}. 
\] 
For prime $p$ and natural number $e$, we will later define a set $\mathrm{U}_p(e,\phi)$ 
which consists of some finite number of complex-valued 
functions on $\Qs_p$ having support in $\Z_p\cap \Qs_p$. 

Let $\Nsc$ be the set of 
positive square-free numbers and 
for $n \in \N$, let $n^{\mathrm{sc}}$ be the square-free part of $n$. 
Let 
$A$ be a function on the set $\Nsc$ having values in $\C$ 
and $E$ be an integer such that
$\widetilde{N_{\phi}} \mid E$. 
Denote $e_p$ = $\nu_p(E)$ for all prime numbers $p$ and 
let $\underline{c} = (c_p)$ be any element 
of $\prod_p{\mathrm{U}_p(e_p,\phi)}$. Define
\[f(\underline{c},A)(z) :=
\sum_{n=1}^{\infty} A(n^{\mathrm{sc}})n^{(k-2)/4}\prod_p{c_p(n)}\ q^n, \qquad z \in \Hup\]
and let $\overline{\mathrm{U}}(E,\phi,A)$ be the space generated by these functions $f(\underline{c},A)$ on 
$\Hup$ where $\underline{c} \in \prod_p{\mathrm{U}_p(e_p,\phi)}$.

With the above notation, we are now ready to state the main theorem of Waldspurger \cite[Th\'{e}or\`{e}me 1]{Waldspurger}.

\begin{thm} \label{thm:wald} (Waldspurger)
Let $\mathrm{(H2)}$ be the following hypothesis:
One of the following holds:
\begin{enumerate}
 \item[(a)] the local component $\rho_2$ is not supercuspidal;
 \item[(b)] the conductor of $\chi_0$ is divisible by $16$;
 \item[(c)] $16 \mid M_{\phi}$.
\end{enumerate}

Let $\chi$ be a Dirichlet character and 
$\phi$ be a newform of weight $k-1$ and character $\chi^2$ such that 
$\mathrm{(H1)}$ and $\mathrm{(H2)}$ hold. Then there exists a function $A_{\phi}$ 
on $\Nsc$ such that for $t \in \Nsc$:
\begin{equation}\label{eq:Aphi}
{A_{\phi}(t)}^{2}:=\mathrm{L}(\phi\otimes \chi_0^{-1} \chi_t, \frac{k-1}{2}) \cdot \epsilon(\chi_0^{-1}\chi_t,\frac{1}{2}).
\end{equation}
Moreover, for $N \geq 1$,
\[
 S_{k/2}(N,\chi,\phi)=\bigoplus \overline{\mathrm{U}}(E,\phi,A_{\phi})
\]
where the sum is 
over all $E \geq 1$ such that $\widetilde{N_{\phi}} \mid E \mid N$.
\end{thm}
Here $\chi_t = \kro{t}{\cdot}$ is a quadratic 
character with conductor $\lvert t \rvert$ if $t \equiv 1 \pmod{4}$, otherwise with conductor 
$\lvert 4t \rvert$ if $t \equiv 2, 3 \pmod{4}$.

\noindent {\bf Remark.} Note that the function $A_\phi$ depends only on $\chi$ and $\phi$. 
However $A_{\phi}$ is not determined by \eqref{eq:Aphi}, so we cannot use this 
theorem for computing a basis for the space $S_{k/2}(N,\chi,\phi)$. 
However, in Theorem~\ref{thm:algo} we have already given an algorithm to
compute this space, and if $f(z) = \sum_{n=1}^{\infty}a_n q^n$ is one of the basis 
elements then we can express the critical value of the $\mathrm{L}$-function of the twist of the newform 
$\phi$ by the character $\chi_0^{-1}\chi_t$, in terms of the square of the Fourier coefficient $a_t$ 
and the factor $\epsilon(\chi_0^{-1}\chi_t,1/2)$ which depends on the local components 
of $\phi$ and $\chi_0$.

It is to be noted that $\epsilon(\chi, 1/2)$ for any Hecke character $\chi$ can be computed 
as shown in Tate's article \cite{Tate}. 
In particular, when $\chi$ is quadratic, $\epsilon(\chi,1/2)$=1. Since we will be dealing only with 
quadratic characters, we can ignore the $\epsilon$-factor. Moreover, note that if $\chi$ is quadratic, then 
the conductor of $\chi_0$ is at most divisible by $8$, so we do not need to 
consider possibility (b) of the hypothesis $\mathrm{(H2)}$.

Further by Corollary \ref{cor:supercuspidal}, possibilities (a) and (c) of the hypothesis $\mathrm{(H2)}$ can be simply 
stated in terms of the level $M_\phi$. Assuming $\chi$ to be quadratic, Waldspurger's Theorem is applicable whenever 
either $M_\phi$ is odd; or $\nu_2(M_\phi)=1$ and $\lambda_2 \ne 0$; or $\nu_2(M_\phi) \geq 4$. The last condition 
is the same as possibility (c) of $\mathrm{(H2)}$.

We also state the following corollary of Waldspurger \cite[p.483]{Waldspurger}.
\begin{cor} (Waldspurger) \label{cor:wald}
Let $\phi \in S_{k-1}^{\mathrm{new}}(M_{\phi},\chi^2)$ be a newform such that $\phi$ satisfies $\mathrm{(H1)}$. 
Suppose~\footnote{
In this corollary we do not require $f$ to be of the form
$f(\underline{c},A_\phi)$.
}
 $f(z) = \sum_{n=1}^{\infty}a_n q^n \in S_{k/2}(N,\chi,\phi)$ 
for some $N \geq 1$ such that $M_{\phi}$ divides $N/2$. Suppose that 
$n_1, n_2 \in \Nsc$ such that $n_1/n_2 \in {\Qs_p}^2$ for all $p \mid N$. Then we have the following 
relation:
\[
a_{n_1}^2 \mathrm{L}(\phi \chi_0^{-1}\chi_{n_2},1)\chi(n_2/n_1)n_2^{{k/2}-1} 
= a_{n_2}^2 \mathrm{L}(\phi \chi_0^{-1}\chi_{n_1},1)n_1^{{k/2}-1}.
\]
\end{cor}

In what follows $( \cdot \ ,\ \cdot )_p$ stands for the Hilbert symbol defined on $\Qs_p \times \Qs_p$. 
Recall that (see for example, \cite{Cohn}) if $p=2$ and $a$, $b$ are odd then
\[(2^s a, 2^t b)_2 = \kro{2}{|a|}^t \kro{2}{|b|}^s (-1)^\frac{(a-1)(b-1)}{4}.\]
For an odd prime $p$ and $a$, $b$ coprime to $p$,
\[(p^s a, p^t b)_p = \kro{-1}{p}^{st} \kro{a}{p}^t \kro{b}{p}^s.\]
In particular, for an odd $n$, $(n,-1)_2 = (-1)^\frac{n-1}{2}$ and $(2,n)_2 = (-1)^\frac{n^2-1}{8}$.
Also, if $\nu_p(n)=0$ then $(p,n)_p = \kro{n}{p}$, and if $\nu_p(n)=1$ and $n = pn^\prime$, 
then $(p,n)_p = \kro{-n^\prime}{p}$.

We now write down explicitly the definitions of the integers $\widetilde{n_p}$ and the
local factors $U(e,\phi)$ used in Waldspurger's Theorem.
It is to be noted that for Waldspurger's Theorem, we require the values of the functions
in $U_p(e,\phi)$ only at square-free positive integers.
We will first define a certain set of functions.
\vskip 3mm

\noindent Case 1. $p\ odd$.\\ 
Waldspurger considered the following set of functions.
\[
\varLambda_p:=\{ c_{p,\delta}^{(0)},\ c_{p,\delta}^{(1)},\ c_{p,\delta}^{(2)},\ c_{p,\delta}^{(3)},\ c_{p,\delta}^{(4)},\ 
c_{p,\delta}^{(5)},\ c_{p,\delta}^{(6)} : \delta \in \C \}.
\]

We will be interested only in values of the functions in $\varLambda_p$ at square-free numbers in $\Z_p\setminus\{0\}$. 
Let $n \in \Z_p\setminus\{0\}$ be square-free, i.e. $\nu_p(n)=0$ or $\nu_p(n)=1$. We get the following 
after simplification:

\[
c_{p,\delta}^{(0)}(n)=1; \]

\begin{eqnarray*}
c_{p,\delta}^{(1)}(n)=\left\{ \begin{array}{ll}
1 & \ \ \text{if}\ \nu_p(n)=0\\
\delta & \ \ \text{if}\ \nu_p(n)=1;\\
\end{array} \right.
\end{eqnarray*}

\begin{eqnarray*}
c_{p,\delta}^{(2)}(n)=\left\{ \begin{array}{ll}
1-(p,n)_p\chi_{0,p}(p)p^{-1/2}\delta^{-1} & \ \ \text{if}\  \nu_p(n)=0\\
1 & \ \ \text{if}\ \nu_p(n)=1;\\
\end{array} \right.
\end{eqnarray*}

\begin{eqnarray*}
c_{p,\delta}^{(3)}(n)=\left\{ \begin{array}{ll}
1 & \ \  \text{if}\ \nu_p(n)=0\\
\delta-(p,n)_p\chi_{0,p}(p)p^{-1/2} & \ \ \text{if}\ \nu_p(n)=1;\\
\end{array} \right.
\end{eqnarray*}

\begin{eqnarray*}
c_{p,\delta}^{(4)}(n)=\left\{ \begin{array}{ll}
0 & \ \ \text{if}\ \nu_p(n)=0\\
\delta(p-1)^{-1} & \ \ \text{if}\ \nu_p(n)=1;\\
\end{array} \right.
\end{eqnarray*}

\begin{eqnarray*}
c_{p,\delta}^{(5)}(n)=\left\{ \begin{array}{ll}
2^{1/2} & \ \ \text{if}\ \nu_p(n)=0 \ \text{and}\  (p,n)_p= -p^{1/2}\chi_{0,p}(p^{-1})\delta\\
0 & \ \ \text{if}\ \nu_p(n)=0 \ \text{and}\ (p,n)_p = p^{1/2}\chi_{0,p}(p^{-1})\delta\\
1 & \ \ \text{if}\ \nu_p(n)=1; \\
\end{array} \right.
\end{eqnarray*}

\begin{eqnarray*}
c_{p,\delta}^{(6)}(n)=\left\{ \begin{array}{ll}
1 & \ \ \text{if}\ \nu_p(n)=0 \\
2^{1/2}\delta & \ \ \text{if}\ \nu_p(n)=1 \ \text{and}\ (p,n)_p= -p^{1/2}\chi_{0,p}(p^{-1})\delta\\
0 & \ \ \text{if}\ \nu_p(n)=1 \ \text{and}\ (p,n)_p = p^{1/2}\chi_{0,p}(p^{-1})\delta.\\
\end{array} \right.
\end{eqnarray*}

\vskip3mm

\noindent Case 2. $p=2$.\\
In this case Waldspurger consider the following set of functions.
\[
\varLambda_2:=\{ c_{2,\delta}^{(0)},\ c_{2,\delta}^{(1)},\ c_{2,\delta}^{(2)},\ c_{2,\delta}^{(3)},\ c_{2,\delta}^{(4)},\ 
c_{2,\delta}^{(5)},\ c_{2,\delta}^{(6)} : \delta \in \C \},
\]
\
Let $n \in \Z_2\setminus\{0\}$ be square-free so that either $\nu_2(n)=0$ or $\nu_2(n)=1$. We have:

\begin{eqnarray*}
c_{2,\delta}^{(0)}(n)=\left\{ \begin{array}{ll}
1 & \ \  \text{if}\ \nu_2(n)=0\\
\delta & \ \  \text{if}\ \nu_2(n)=1; \\
\end{array} \right.
\end{eqnarray*}

\begin{eqnarray*}
c_{2,\delta}^{(1)}(n)=\left\{ \begin{array}{ll}
\delta-(2,n)_2\chi_{0,2}(2)2^{-1/2} & \ \  \text{if}\ \nu_2(n)=0\ \text{and}\ (n,-1)_2= \chi_{0,2}(-1) \\
1 & \ \ \text{if}\ \nu_2(n)=0\ \text{and}\ (n,-1)_2 = -\chi_{0,2}(-1) \\
1 & \ \ \text{if}\ \nu_2(n)=1; \\
\end{array} \right.
\end{eqnarray*}

\begin{eqnarray*}
c_{2,\delta}^{(2)}(n)=\left\{ \begin{array}{ll}
\delta & \ \ \text{if}\ \nu_2(n)=0\ \text{and}\ (n,-1)_2= \chi_{0,2}(-1)\\
0 & \ \ \text{if}\ \nu_2(n)=0\ \text{and}\ (n,-1)_2 = -\chi_{0,2}(-1)\\
0 & \ \ \text{if}\ \nu_2(n)=1; \\
\end{array} \right.
\end{eqnarray*}

\begin{eqnarray*}
c_{2,\delta}^{(3)}(n)=\left\{ \begin{array}{ll}
\delta^{-1} & \ \ \text{if}\ \nu_2(n)=0\\
\delta-(2,n)_2\chi_{0,2}(2)2^{-1/2} & \ \ \text{if}\ \nu_2(n)=1\ \text{and}\ (n,-1)_2= \chi_{0,2}(-1) \\
1 & \ \ \text{if}\ \nu_2(n)=1\ \text{and}\ (n,-1)_2 = -\chi_{0,2}(-1); \\
\end{array} \right.
\end{eqnarray*}

\begin{eqnarray*}
c_{2,\delta}^{(4)}(n)=\left\{ \begin{array}{ll}
0 & \ \ \text{if}\ \nu_2(n)=0\\
2 \delta-(2,n)_2\chi_{0,2}(2)2^{-1/2} & \ \ \text{if}\ \nu_2(n)=1\ \text{and}\ (n,-1)_2= \chi_{0,2}(-1) \\
1 & \ \ \text{if}\ \nu_2(n)=1\ \text{and}\ (n,-1)_2 = -\chi_{0,2}(-1);\\
\end{array} \right.
\end{eqnarray*}

\begin{eqnarray*}
c_{2,\delta}^{(5)}(n)=\left\{ \begin{array}{ll}
0 & \ \ \text{if}\ \nu_2(n)=0,\ (n,-1)_2= \chi_{0,2}(-1)\ \text{and}\ (2,n)_2= 2^{1/2}\chi_{0,2}(2^{-1})\delta\\
2^{1/2}\delta & \ \ \text{if}\ \nu_2(n)=0,\ (n,-1)_2= \chi_{0,2}(-1)\ \text{and}\ (2,n)_2 = -2^{1/2}\chi_{0,2}(2^{-1})\delta\\
1 & \ \ \text{if}\ \nu_2(n)=0\ \text{and}\ (n,-1)_2 = -\chi_{0,2}(-1)\\
1 & \ \ \text{if}\ \nu_2(n)=1; \\
\end{array} \right.
\end{eqnarray*}

\begin{eqnarray*}
c_{2,\delta}^{(6)}(n)=\left\{ \begin{array}{ll}
\delta^{-1} & \ \ \text{if}\ \nu_2(n)=0 \\
0 & \ \ \text{if}\ \nu_2(n)=1,\ (n,-1)_2= \chi_{0,2}(-1)\ \text{and}\ (2,n)_2 = 2^{1/2}\chi_{0,2}(2^{-1})\delta\\
2^{1/2}\delta & \ \ \text{if}\ \nu_2(n)=1,\ (n,-1)_2= \chi_{0,2}(-1)\ \text{and}\ (2,n)_2 = -2^{1/2}\chi_{0,2}(2^{-1})\delta\\
1 & \ \ \text{if}\ \nu_2(n)=1\ \text{and}\ (n,-1)_2= -\chi_{0,2}(-1).\\
\end{array} \right.
\end{eqnarray*}

We will be interested in the above functions only for particular values of $\delta$. 
We will specify and further simplify them later.

Recall that $\lambda_p$ is the Hecke eigenvalue of $\phi$ corresponding
to the Hecke operator $T_p$ for any prime $p$, and $m_p = \nu_p(M_\phi)$. Let $\lambda_p^{\prime} = p^{1-k/2}\lambda_p$.
For $p\nmid M_{\phi}$ let $\alpha_p$ and $\alpha_p^{\prime}$ be such that
\[
 \alpha_p + \alpha_p^{\prime} = \lambda_p^{\prime},
\]
\[
 \alpha_p\cdot\alpha_p^{\prime} = 1.
\]
It is to be noted that if $\phi$ is rational newform of weight $2$ then $\alpha_p \ne \alpha_p^{\prime}$, 
since otherwise $\lambda_p^2 = 4$, which is a contradiction as $\lambda_p$ is rational ($p$-th 
Fourier coefficient of $\phi$).

Next, we need to consider a subset of ${\Qs_p}/{\Qs_p}^2$, denoted by $\Omega_p(\phi)$, which is defined as \label{omega}
\begin{equation}\label{eq:Omega_p}
\begin{split}
\Omega_p(\phi)=  \{& \omega \in {\Qs_p}/{\Qs_p}^2 : \ \exists\ f \in S_{k/2}(N,\chi,\phi)\ 
\text{ for some $N$ and $\exists \ n \geq 1$ such} \\ 
& \text{that }
 i)\ \text{the image of $n$ in}\ {\Qs_p}/{\Qs_p}^2\ \text{is}\ \omega\ ;\ ii)\ \text{the $n$th}\ \text{coefficient of }f \ne 0\}.
\end{split}
\end{equation}
Note that the set $\Omega_p(\phi)$ depends on the newform $\phi$ and character $\chi$ that we started with. 
Computation of this set is important in our applications and we will see that we need this set only in the case
when $m_p \geq 1$ and $\lambda_p=0$. Since this set consists of at most eight elements when $p=2$,  
and four when $p$ is an odd prime, computation doesn't seem to be difficult. Indeed, we can use the results of 
Section~\ref{section:coeffmodulon} and our algorithm in Theorem~\ref{thm:algo} to compute most of the elements. 

Waldspurger defined another set of local functions on ${\Qs_p}/{\Qs_p}^2$ taking values in $\Z/2\Z$:
\[
\Gamma_p:=\{ \gamma_{e,\upsilon} : e \in \Z,\ \upsilon \in {\Qs_p}/{\Qs_p}^2 \mathrm{such\ that\ } 
\nu_p(\upsilon)\equiv e\ (\mathrm{mod\ }2)\},
\]
where 
\begin{eqnarray*}
\gamma_{e,\upsilon}(u)=\left\{ \begin{array}{ll}
1 & \ \ \text{if}\ u \in \upsilon{\Qs_p}^2\ \text{and}\ \nu_p(u)=e\\
0 & \ \ \mathrm{else}.\\
\end{array} \ \ \right.
\end{eqnarray*}

If $p=2$, define 
\[
\gamma_{e,\upsilon}^{\prime} = \frac{1}{2}(\gamma_{e,\upsilon} + \gamma_{e,5\upsilon}),
 \]

\begin{eqnarray*}
\gamma_{e}^{''}(u)=\left\{ \begin{array}{ll}
1 & \ \ \text{if}\ \nu_2(u)=e\\
0 & \ \ \mathrm{else},\\
\end{array} \ \ \right.
\end{eqnarray*}
and 
\begin{eqnarray*}
\gamma_{e}^{0}(u)=\left\{ \begin{array}{ll}
1 & \ \ \text{if}\ \nu_2(u)=e\ \text{and}\ (u,-1)_2 = -\chi_{0,2}(-1)\ \mathrm{or}\ \nu_2(u)=e+1 \\
0 & \ \ \mathrm{else}.\\
\end{array} \ \ \right.
\end{eqnarray*}

Now we are ready to define the local factors $\widetilde{n_p}$ and the set $U_p(e,\phi)$ for
$e = \widetilde{n_p}$. We will be dealing with several cases and subcases and in each of them
we will be simplifying Waldspurger's formulae and making them more explicit for our use.

\begin{enumerate}[leftmargin=1.5cm]
\item[Case 1.] $p\ odd$ and $m_p \geq 1$. 

We consider the following subcases:
\begin{enumerate}[leftmargin=-0.25cm]
\item[(a)] $\lambda_p = 0$.

In this case we need to compute $\Omega_p(\phi)$. We know that ${\Qs_p}/{\Qs_p}^2 =\{ 1,\ p,\ u,$ $\ pu \}$ 
where $u$ is unit in $\Z_p$ which is a non-square mod $p$. 
If there exists a $\omega \in \Omega_p(\phi)$ such that $\nu_p(\omega) = 0$ then $\widetilde{n_p} = m_p$, and
$U_p(\widetilde{n_p},\phi) = \{\gamma_{0,\omega} :\ \omega \in \Omega_p(\phi)\ \mathrm{and}\ \nu_p(\omega) = 0\}$.
In this case, the set $U_p(\widetilde{n_p},\phi)$ consists of at most 
the functions $\gamma_{0,1}$ and $\gamma_{0,u}$.
Otherwise, for all $\omega \in \Omega_p(\phi)$, $\nu_p(\omega) = 1$. In this case $\widetilde{n_p} = m_p+1$, and
$U_p(\widetilde{n_p},\phi) = \{\gamma_{1,\omega} :\ \omega \in \Omega_p(\phi)\ \mathrm{and}\ \nu_p(\omega) = 1\}$, hence 
$U_p(\widetilde{n_p},\phi)$ consists of at most $\gamma_{1,p}$ and $\gamma_{1,pu}$.
Note that $\gamma_{0,1},\ \gamma_{0,u},\ \gamma_{1,p},\ \gamma_{1,pu}$ are 
characteristic functions of $1,\ u,\ p,\ pu$ modulo ${\Qs_p}^2$ respectively.
\item[(b)] $\lambda_p \ne 0$.

In this case we must have $m_p =1 $, since $m_p \geq 2$ implies 
that $\lambda_p = 0$. Note that $p \in S$ since by Corollary \ref{cor:special} 
$\rho_p$ is a special representation and hence not irreducible principal series.
We have further subcases:

\begin{enumerate}[leftmargin=0.25cm]
\item[(i)] $\chi_{0,p}$ is unramified.

Here again $\widetilde{n_p} = m_p =1$ and $U_p(1,\phi)=\{c_{p,\lambda_p^{\prime}}^{(5)}\}$.
We use the theory of newforms to simplify the function $c_{p,\lambda_p^{\prime}}^{(5)}$. 
Since $m_p =1$ we get that $\lambda_p = -\omega_p p^{(k-3)/2}$ and $\lambda_p^{\prime} = -\omega_p p^{-1/2}$. 
Here $\omega_p \in \{\pm1\}$ is the eigenvalue under the Atkin-Lehner involution corresponding to the prime $p$.
Hence we have in this case,
\begin{eqnarray*}
c_{p,\lambda_p^{\prime}}^{(5)}(n)=\left\{ \begin{array}{ll}
2^{1/2} & \ \ \text{if}\ \nu_p(n)=0\ \text{and}\ \kro{n}{p}= \omega_p\chi_{0,p}(p^{-1})\\
0 & \ \ \text{if}\ \nu_p(n)=0\ \text{and}\ \kro{n}{p} = -\omega_p\chi_{0,p}(p^{-1})\\
1 & \ \ \text{if}\ \nu_p(n)=1. \\
\end{array} \right.
\end{eqnarray*}
\item[(ii)]$\chi_{0,p}$ is ramified.

We have $\widetilde{n_p} = m_p =1$ and $U_p(1,\phi)=\{c_{p,\lambda_p^{\prime}}^{(6)}\}$.
As in the above subcase, we get the following simplification:
\begin{eqnarray*}
c_{p,\lambda_p^{\prime}}^{(6)}(n)=\left\{ \begin{array}{ll}
1 & \ \ \text{if}\ \nu_p(n)=0 \\
-\omega_p 2^{1/2}p^{-1/2} & \ \ \text{if}\ \nu_p(n)=1\ \text{and}\ (p,n)_p= \omega_p \chi_{0,p}(p^{-1})\\
0 & \ \ \text{if}\ \nu_p(n)=1\ \text{and}\ (p,n)_p = -\omega_p \chi_{0,p}(p^{-1}).\\
\end{array} \right.
\end{eqnarray*}
\end{enumerate}
\end{enumerate}
\item[Case 2.] $p\ odd$ and $m_p =0$. 

We have the following subcases:
\begin{enumerate}[leftmargin=-0.25cm]
\item[(a)] $\chi_{0,p}$ is unramified.

Here, $\widetilde{n_p} = m_p =0$ and $U_p(0,\phi)=\{c_{p,\lambda_p^{\prime}}^{(0)}\}$.
It is to be noted that $c_{p,\lambda_p^{\prime}}^{(0)}$ takes the value 1 at any square-free $n$.

\item[(b)] $\chi_{0,p}$ is ramified.

We have $\widetilde{n_p} = 1$ and $U_p(1,\phi)=\{c_{p,\alpha_p}^{(3)},\ c_{p,\alpha_p^{\prime}}^{(3)}\}$ if 
$\alpha_p \ne \alpha_p^{\prime}$, else $U_p(1,\phi)=\{c_{p,\alpha_p}^{(3)},\ c_{p,\alpha_p}^{(4)}\}.$
We note that if $p$ does not divide the modulus of $\chi$, then we do not need to consider this subcase
because in this case $\chi_{0,p}$ is unramified by Lemma \ref{lem:charprop}.
\end{enumerate}

\item[Case 3.] $p=2$ and $m_2 \geq 1$. 

Consider the following subcases:

\begin{enumerate}[leftmargin=-0.25cm]
\item[(a)] $\lambda_2 = 0$.

We compute $\Omega_2(\phi)$. Note that ${\Qs_2}/{\Qs_2}^2 =\{ \pm1,\ \pm2,\ \pm5,\ \pm10 \}$.
If there exists a $\omega \in \Omega_2(\phi)$ such that $\nu_2(\omega) = 0$ then $\widetilde{n_2} = m_2+2$, and
$U_2(\widetilde{n_2},\phi) = \{\gamma_{0,\omega} :\ \omega \in \Omega_2(\phi)\ \mathrm{and}\ \nu_2(\omega) = 0\}$.
In this case, the set $U_2(\widetilde{n_2},\phi)$ consists of at most $\gamma_{0,1}$, $\gamma_{0,3}$, $\gamma_{0,5}$, 
and $\gamma_{0,7}$ .
Otherwise, for all $\omega \in \Omega_2(\phi)$, $\nu_2(\omega) = 1$ and then $\widetilde{n_2} = m_2+3$, and
$U_2(\widetilde{n_2},\phi) = \{\gamma_{1,\omega} :\ \omega \in \Omega_2(\phi)\ \mathrm{and}\ \nu_2(\omega) = 1\}$, hence 
$U_2(\widetilde{n_2},\phi)$ consists of at most $\gamma_{1,2}$, $\gamma_{1,6}$, $\gamma_{1,10}$ and $\gamma_{1,14}$.
As above, $\gamma_{0,i}$ for $i \in \{1,\ 3,\ 5,\ 7\}$ are the characteristic functions of an odd residue class modulo $8$ and 
$\gamma_{1,j}$ for $j \in \{2,\ 6,\ 10,\ 14\}$ are the characteristic functions of even residue class modulo ${\Qs_2}^2$.

\item[(b)] $\lambda_2 \ne 0$.

We must have $m_2 =1 $. As before Corollary \ref{cor:special} implies that $\rho_2$ is a special 
representation and hence $p \in S$. We have the following subcases:
\begin{enumerate}[leftmargin=0.25cm]
\item[(i)] $\chi_{0,2}$ is trivial on $1 + 4\Z_2$.

Here $\widetilde{n_2} = 2$ and $U_2(2,\phi)=\{c_{2,\lambda_2^{\prime}}^{(5)}\}$.
Since $m_2 =1$ we get that $\lambda_2 = -\omega_2 2^{(k-3)/2}$ and $\lambda_2^{\prime} = -\omega_2 2^{-1/2}$; 
$\omega_2 \in \{ \pm1\}$ is the eigenvalue under the Atkin-Lehner involution corresponding to $2$.
Hence we have,
\begin{eqnarray*}
c_{2,\lambda_2^{\prime}}^{(5)}(n)=\left\{ \begin{array}{ll}
0 & \ \ \text{if}\ \nu_2(n)=0 ,\ (-1)^{\frac{n-1}{2}}= \chi_{0,2}(-1)\ \text{and}\ (-1)^{\frac{n^2-1}{8}}= -\omega_2 \chi_{0,2}(2^{-1})\\
-\omega_2  & \ \ \text{if}\ \nu_2(n)=0,\ (-1)^{\frac{n-1}{2}}= \chi_{0,2}(-1)\ \text{and}\ (-1)^{\frac{n^2-1}{8}} = \omega_2 \chi_{0,2}(2^{-1})\\
1 & \ \ \text{if}\ \nu_2(n)=0,\ (-1)^{\frac{n-1}{2}} = -\chi_{0,2}(-1)\\
1 & \ \ \text{if}\ \nu_2(n)=1. \\
\end{array} \right.
\end{eqnarray*}

\item[(ii)]$\chi_{0,2}$ is nontrivial on $1 + 4\Z_2$.

Here $\widetilde{n_2} = 3$ and $U_2(3,\phi)=\{c_{p,\lambda_2^{\prime}}^{(6)},\ \gamma_{0}^{''}\}$ and 
we get the following simplification:
\begin{eqnarray*}
c_{2,\lambda_2^{\prime}}^{(6)}(n)=\left\{ \begin{array}{ll}
-\omega_2 2^{1/2} & \ \ \text{if}\ \nu_2(n)=0 \\
0 & \ \ \text{if}\ \nu_2(n)=1,\ (n,-1)_2= \chi_{0,2}(-1)\ \text{and}\ (2,n)_2 = -\omega_2 \chi_{0,2}(2^{-1})\\
-\omega_2 & \ \ \text{if}\ \nu_2(n)=1,\ (n,-1)_2= \chi_{0,2}(-1)\ \text{and}\ (2,n)_2 = \omega_2 \chi_{0,2}(2^{-1})\\
1 & \ \ \text{if}\ \nu_2(n)=1,\ (n,-1)_2= -\chi_{0,2}(-1).\\
\end{array} \right.
\end{eqnarray*}
\end{enumerate}
\end{enumerate}

\item[Case 4.] $p=2$ and $m_2 =0$. 

We have the following subcases:
\begin{enumerate}[leftmargin=-0.25cm]
\item[(a)] $\chi_{0,2}$ is trivial on $1 + 4\Z_2$.

We have $\widetilde{n_2} = 2$ and $U_2(2,\phi)=\{c_{2,\alpha_2}^{(1)},\ c_{2,\alpha_2^{\prime}}^{(1)}\}$
if $\alpha_2 \ne \alpha_2^{\prime}$, else $U_2(2,\phi)=\{c_{2,\alpha_2}^{(1)},\ c_{2,\alpha_2}^{(2)}\}$.

\item[(b)] $\chi_{0,2}$ is nontrivial on $1 + 4\Z_2$.

Here $\widetilde{n_2} = 3$ and $U_2(3,\phi)=\{c_{2,\alpha_2}^{(3)},\ c_{2,\alpha_2^{\prime}}^{(3)},\ \gamma_{0}^{''}\}$ if 
$\alpha_2 \ne \alpha_2^{\prime}$, else $U_2(3,\phi)=\{c_{2,\alpha_2}^{(3)},\ c_{2,\alpha_2}^{(4)},\ \gamma_{0}^{''}\}.$
\end{enumerate}
\end{enumerate}

We would like to point out the following useful lemma:
\begin{lem}
Let $\chi$ be a quadratic character modulo $N$ such that $\nu_2(N)$ is at most $2$. Then,
$\chi_{0,2}$ is trivial on $1 + 4\Z_2$.
\end{lem}
\begin{proof}
Since $\chi$ is a quadratic, $\chi_0$ is also quadratic with modulus 
$\LCM(4,N) = 4N^\prime$ where $2 \nmid N^\prime$. Now the
lemma follows from part (iii) of Lemma~\ref{lem:charprop}.
\end{proof}

\noindent {\bf Remark.}
These simplifications along with our method to compute a basis for $S_{k/2}(N,\chi,\phi)$ for 
suitable $N$ and $\chi$ lead to an algorithm for computing 
critical values of the $\mathrm{L}$-functions of 
certain quadratic twists of $\phi$. For example, if $M_\phi = p^\alpha$ for some odd prime $p$, then the 
possible values for $\widetilde{N_{\phi}}$ are
 either $4p^\alpha$ or $4p^{\alpha+1}$, hence we compute 
bases for spaces $S_{k/2}(4p^\alpha, \chi_{\mathrm{triv}}, \phi)$ and 
$S_{k/2}(4p^{\alpha+1},\chi_{\mathrm{triv}},\phi)$ 
and the sets $U_2(2,\phi)$, $U_p(\alpha,\phi)$, $U_p(\alpha+1,\phi)$ to apply 
Theorem \ref{thm:wald} in order to get the desired results.

Note that in the above we have discussed computation of $U_p(e,\phi)$ only 
for $e = \widetilde{n_p}$. But in certain cases working with the level $\widetilde{N_{\phi}}$ is not sufficient to get 
the complete information and one might need to go to higher levels.

\section{Periods}\label{sec:period}

\begin{lem}\label{lem:period}
Let $E/\Q$ be an elliptic curve, given by a minimal Weierstrass model,
and let $E_n$ be the minimal model of its twist by a square-free
positive integer $n$. Then there is a computable non-zero rational number
$\alpha_n$ such that
\[
\Omega(E_n)=\frac{\alpha_n \Omega(E)}{\sqrt{n}}.
\]
\end{lem}
The proof we give also explains how to compute $\alpha_n$.
\begin{proof}

Let $\omega=dx/(2y+a_1 x +a_3)$ be the invariant differential for the model 
\[
E: y^2+a_1 xy +  a_3 y= x^3+a_2 x^2+a_4 x +a_6.
\]
By definition, the period
\[
\Omega(E)= \int_{E(\R)} \lvert \omega \rvert.
\]
Recall \cite[p.49]{Silverman} that a change of variable
\[
x=u^2 x^\prime +r, \qquad y= u^3 y^\prime+ u^2 s x^\prime+t
\]
leads to a model $E^\prime$ with invariant differential $\omega^\prime= u \omega$;
thus the periods are related by $\Omega(E^\prime)= \lvert u \rvert \Omega(E)$.
Completing the square in $y$ we obtain the model
\[
E^\prime : {y^\prime}^2= {x^\prime}^3+A {x^\prime}^2 + B {x^\prime} + C
\]
where
\[
A=\frac{b_2}{4}, \qquad B= \frac{b_4}{2}, \qquad C=\frac{b_6}{4}.
\]
Since $u=1$ in this change of variable, $\omega^\prime=\omega$
and $\Omega(E^\prime)=\Omega(E)$.
Now let the model $E^{\prime\prime}$ be the twist of $E^\prime$ by $n$:
\[
E^{\prime\prime} : {y^{\prime\prime}}^2= {x^{\prime\prime}}^3+A n {x^{\prime\prime}}^2 + B n^2 {x^{\prime\prime}} + Cn^3.
\]
Note that these are related by the change of variable
\[
y^{\prime\prime} = n^{3/2}  y^\prime, \qquad x^{\prime\prime}= n x^\prime.
\]
Thus the invariant differentials satisfy
\[
\omega^{\prime\prime} = \frac{d x^{\prime\prime}}{2 y^{\prime\prime}} = \frac{\omega^\prime}{\sqrt{n}}.
\]
Thus
\[
\Omega(E^{\prime\prime})=\frac{\Omega(E^\prime)}{\sqrt{n}}=\frac{\Omega(E)}{\sqrt{n}}.
\]
Now the model $E^{\prime\prime}$ is not necessarily minimal (nor even integral at $2$), but by Tate's
algorithm there is a change
of variables
\[
x^{\prime\prime}=u^2 X +r, \qquad y^{\prime\prime}= u^3 Y+ u^2 s X+t
\]
with rational $u$, $s$, $t$ (and $u \neq 0$) such that the resulting model $E_n$ is minimal.
By the above
\[
\Omega(E_n)=u \Omega(E^{\prime\prime})= \frac{\lvert u \rvert \Omega(E)}{\sqrt{n}}.
\]
\end{proof}

\begin{lem}\label{lem:celery}
Let $E : Y^2=X^3+AX^2+BX+C$ be an elliptic curve 
with $A$, $B$, $C \in \Z$. Suppose that the discriminant
of this model is sixth-power free. Let
$n$ be a square-free positive integer. Then a minimal model
for the $n$-th twist is $E_n : Y^2=X^3+An X^2+B n^2 X + Cn^3$.
Moreover, the periods are related by the formula
\[
\Omega(E_{n})=\frac{\Omega(E_{1})}{\sqrt{n}}.
\]
\end{lem}
\begin{proof}
Let $\Delta$ be the discriminant of the model $E : Y^2=X^3+AX^2+BX+C$.
We are assuming that $\Delta$ is sixth-power free. Thus it is $12$-th power
free, and so $E$ is minimal. Now the model
$E_n : Y^2=X^3+An X^2+B n^2 X + Cn^3$ has discriminant
$\Delta_n=\Delta \cdot n^6$. Since $n$ is square-free this is $12$-th power free.
Thus the model for $E_n$ is minimal. The argument in the proof of
Lemma~\ref{lem:period} completes the proof.
\end{proof}

\section{Modular Forms are Determined by Coefficients Modulo n}\label{section:coeffmodulon}

As usual $N$ is a positive integer divisible by $4$, $\chi$ a Dirichlet
character modulo $N$. Let $k$ be an odd integer. 
Let $\phi$ be a newform of weight $k-1$, level dividing $N/2$ and
character $\chi^2$. To apply Waldspurger's Theorem,
we need to know (see p.\pageref{omega}) for certain primes $p$,  certain
$\omega \in {\Qs_p}/{\Qs_p}^2$ and certain 
forms $f=\sum a_n q^n \in S_{k/2}(N,\chi,\phi)$, whether there is 
some $n$ such that the image of $n$ in ${\Qs_p}/{\Qs_p}^2$
is $\omega$ and $a_n \ne 0$. Given such $p$, $f$ and $\omega$
we can write down the first few coefficients of $f$ and test
whether the image of $n$ in ${\Qs_p}/{\Qs_p}^2$ is $\omega$ and
$a_n \ne 0$. If there is such an $n$ then we should be able
to find it by writing down enough coefficients. However,
sometimes it appears that $a_n=0$ for all $n$
that are equivalent in ${\Qs_p}/{\Qs_p}^2$ to $\omega$.
To be able to prove that, we have developed the results in  this section.

\begin{thm}\label{thm:congruence}
Let $N$ be a positive integer such that $4\mid N$ and $\chi$ be a Dirichlet 
character modulo $N$. Let $f(z) = \sum_{n=1}^{\infty} a_n q^n \in S_{k/2}(N,\chi)$.
Let $a$, $M$ be integers such that $(a,M)=1$. Let $R = \frac{k}{24} [ \SL_2(\Z) : \mathrm{\Gamma_1}(NM^2) ]$.
Suppose $a_n = 0$ whenever $n \not\equiv a \pmod{M}$ for all integers $n$ up to $R+1$. Then
$a_n = 0$ whenever $n \not\equiv a \pmod{M}$ for all $n$. 
Moreover if $M^2 \mid N$ then the above statement holds with
\[R = 
\begin{cases} \frac{k}{24} [ \SL_2(\Z) : \mathrm{\Gamma_1}(N) ] & \text{if $\frac{N}{M} \equiv 0 \pmod{4}$}\\
 \frac{k}{24} [ \SL_2(\Z) : \mathrm{\Gamma_1}(2N) ] & \text{if $\frac{N}{M} \equiv 2 \pmod{4}$}.
\end{cases}
\]
\end{thm}

We will be requiring the analogue, in the case of half-integral weight forms, 
of the following theorem of Sturm.
\begin{thm} (Sturm~\cite[p.276]{Sturm}) \label{thm:Sturm}
Let $\mathrm{\Gamma}$ be a congruence subgroup and $k$ be a positive integer. 
Let $f$, $g \in M_k(\mathrm{\Gamma})$ such that $f$ and $g$ have coefficients in $\OO_F$, 
the ring of integers of a number field $F$. Let $\lambda$ be a prime ideal of $\OO_F$.
If 
\[
 \ord_{\lambda}(f - g) > \frac{k}{12} [ \SL_2(\Z) : \mathrm{\Gamma} ],
\]
then $\ord_{\lambda}(f - g) = \infty$, i.e., $f \equiv g \pmod{\lambda}$.
\end{thm}
In the above statement if $f(z) = \sum_{n\geq 0} a_nq^n$ 
then $\ord_{\lambda}(f) := \mathrm{inf}\{n : a_n \notin \lambda \}$. 
If $a_n \in \lambda$ for all $n$, then we let $\ord_{\lambda}(f) := \infty$.

\begin{lem}\label{lem:Sturm}
Let $\mathrm{\Gamma^{\prime}}$ be a congruence subgroup such that $\mathrm{\Gamma^{\prime}} 
\subseteq \mathrm{\Gamma_0}(4)$ and $k^{\prime}$ be a positive odd integer. Then the statement
of Theorem~\ref{thm:Sturm} is valid for $\mathrm{\Gamma} = \mathrm{\Gamma^{\prime}}$ and
$k = k^{\prime}/2$.
\end{lem}

\begin{proof}
Let $h := f - g \in S_{k^{\prime}/2}(\mathrm{\Gamma^{\prime}})$. By assumption, 
$\ord_{\lambda}(h) > \frac{k^{\prime}}{24} [ \SL_2(\Z) : \mathrm{\Gamma^{\prime}} ]$.
Let $h^{\prime}=h^4$. Then $h^{\prime} \in M_{2k^{\prime}}(\mathrm{\Gamma^{\prime}})$. This is because for any 
$\gamma = \left[ \begin{smallmatrix} a&b\\ c&d \end{smallmatrix} \right] \in \Gamma^{\prime}$
and $z \in \Hup$,
\begin{equation*}
\begin{split}
h^{\prime}(\gamma z) &= h^4(\gamma z) \\
&=j(\gamma,z)^{4k^{\prime}} h^4(z)\\
&=(cz+d)^{2k^{\prime}} h^{\prime}(z). \\
\end{split}
\end{equation*}
Also, $\ord_{\lambda}(h^{\prime}) = 4\cdot \ord_{\lambda}(h) > \frac{2k^{\prime}}{12} 
[ \SL_2(\Z) : \mathrm{\Gamma^{\prime}} ]$. So we apply Theorem~\ref{thm:Sturm} to $h^{\prime}$ 
to get that $\ord_{\lambda}(h^{\prime})=\infty$. Hence $\ord_{\lambda}(h) =\infty$. 
\end{proof}
We note that the above lemma still holds if $f$, $g \in M_{k^\prime/2}(\Gamma_0(N),\chi)$; the above 
proof works by taking $h^{\prime}=h^{4n}$ where $n$ is the order of Dirichlet character $\chi$.

We will need the following lemmas for the proof of Theorem~\ref{thm:congruence}.
\begin{lem}
Let $M$ be a positive integer and $a \in \Z$ such that $(a,M)=1$. Define
\[
\mathrm{I_a}(n):= 
\begin{cases}
1  & \text{if $n \equiv a \pmod{M}$}\\
0  & \text{otherwise}.
\end{cases}
\]
Then we have 
\[
 \mathrm{I_a}(n) = \sum_{\psi \in \mathrm{X}(M)} \frac{\psi(a)^{-1}}{\varphi(M)} \psi(n)
\]
where $\mathrm{X}(M)$ denotes the group of Dirichlet characters of modulus $M$ and $\varphi$
is Euler's phi function.
\end{lem}
\begin{proof}
See \cite[p.63, Chapter 6]{Serre}.
\end{proof}

\begin{lem}\label{lem:con}
Let $\left[ \begin{smallmatrix} a&b\\ c&d \end{smallmatrix} \right] \in \Gamma_0(N)$ and $m^2 \mid N$.
Let $0\leq \nu^{\prime} < m$ and $\frac{c\nu^{\prime}}{m} \equiv 0 \pmod{4}$. Then,
$\kro{c}{d+c\nu^{\prime}/m} = \kro{c}{d}$.
\end{lem}
The proof of the above lemma requires the following reciprocity law as stated in   
Cassels and Fr\"{o}hlich \cite[p.350]{CF}: 
\begin{prop}\label{prop:reciprocity}
Let $P$, $Q$ be positive odd integers and $a$ be any non-zero integer with $a= 2^{\alpha}a_0$, $a_0$ odd. Then,
\[\kro{a}{P} = \kro{a}{Q} \text{if $P \equiv Q \pmod{8 a_0}$}.\]
\end{prop}

\begin{prop}\label{prop:twisthalf}
Let $k$ be a positive odd integer, $\chi$ be a Dirichlet character modulo $N$ where $4\mid N$ and 
$f(z) = \sum_{n=0}^{\infty} a_n q^n \in M_{k/2}(N,\chi)$. Suppose $\psi$ is a Dirichlet character 
of conductor $m$ and $f_{\psi}(z) = \sum_{n=0}^{\infty} \psi(n)a_n q^n$. Then,
\begin{itemize}
\item[(i)] $f_\psi \in  M_{k/2}(Nm^2,\chi{\psi}^2)$.
\item[(ii)] If $m^2 \mid N$ and $\frac{N}{m} \equiv 0 \pmod{4}$ then $f_\psi \in M_{k/2}(N,\chi{\psi}^2)$.
\item[(iii)] If $m^2 \mid N$ and $\frac{N}{m} \equiv 2 \pmod{4}$ then $f_\psi \in M_{k/2}(2N,\chi{\psi}^2)$.
\end{itemize}
Moreover, if $f$ is a cusp form then so is $f_{\psi}$.
\end{prop}
\begin{proof}
The proof essentially follows that of Proposition 17 in \cite[Chapter III]{Kob},
which is the integral weight case, with some necessary changes. We use Lemma~\ref{lem:con} to obtain (ii) and (iii).
\end{proof}
\begin{lem}\label{lem:con1}
Let $k$, $N$ be positive integers such that $k$ is odd and $4 \mid N$. Suppose 
$f(z) = \sum_{n=1}^{\infty} a_n q^n \in S_{k/2}(N,\chi)$. Let $a$, $M$ 
be positive integers such that $(a,M)=1$. Define
\[
g(z) := \sum_{n=1}^{\infty} \mathrm{I_a}(n) a_nq^n .
\]
Then $g \in S_{k/2}(\mathrm{\Gamma_1}(NM^2))$.
\end{lem} 
\begin{proof}
We have
\begin{equation*}
\begin{split}
g(z) &= \sum_{n=1}^{\infty} \mathrm{I_a}(n) a_nq^n\\
&= \sum_{n=1}^{\infty} \sum_{\psi \in \mathrm{X}(M)} \frac{\psi(a)^{-1}}{\varphi(M)} \psi(n)\ a_nq^n\\
&= \sum_{\psi \in \mathrm{X}(M)} \alpha_{\psi} \sum_{n=1}^{\infty} \psi(n)\ a_nq^n\\
&= \sum_{\psi \in \mathrm{X}(M)} \alpha_{\psi} f_{\psi} \ ,
\end{split}
\end{equation*}
where $\alpha_{\psi} = \frac{\psi(a)^{-1}}{\varphi(M)}$. Using Proposition~\ref{prop:twisthalf}, 
for all $\psi \in \mathrm{X}(M)$ we have $f_{\psi} \in S_{k/2}(\mathrm{\Gamma_1}(NM^2))$. 
Hence $g \in S_{k/2}(\mathrm{\Gamma_1}(NM^2))$.
\end{proof}

Now we are ready to prove Theorem \ref{thm:congruence}.

\begin{proof}[Proof of Theorem~\ref{thm:congruence}]

Let $h = f-g$ where $g$ is as in the above lemma. Since
$f \in S_{k/2}(\mathrm{\Gamma_1}(NM^2))$, so does $h$. It is clear that
\[\text{coefficient of $q^{n}$ in $h$} = 
\begin{cases}
a_n & \text{if $n \not\equiv a \pmod{M}$}\\
0 & \text{otherwise}.
\end{cases} 
 \]
Thus, $h(z) = \sum \limits_ {n \not\equiv a \pmod{M}} a_nq^n \in S_{k/2}(\mathrm{\Gamma_1}(NM^2))$.
Since we have assumed $a_n = 0$ whenever $n \not\equiv a \pmod{M}$ for all integers $n$ 
up to $R+1$, we apply Lemma~\ref{lem:Sturm} to get $h=0$. If $M^2 \mid N$ we apply parts (ii) and (iii) of
Proposition~\ref{prop:twisthalf} to Lemma \ref{lem:con1}.
\end{proof}

\noindent {\bf Remark.} Note that in Lemma~\ref{lem:con1} if all the Dirichlet characters modulo $M$ 
are quadratic then by Proposition~\ref{prop:twisthalf}, in fact, $g \in S_{k/2}(\mathrm{\Gamma_0}(NM^2, \chi))$. 
Hence in this case Theorem~\ref{thm:congruence} holds with $R = \frac{k}{24} [ \SL_2(\Z) : \mathrm{\Gamma_0}(NM^2) ]$.
For example, if $N=1984$, $k=3$ and $M=8$, since all Dirichlet characters modulo $8$ are quadratic we have $R = 
\frac{3}{24} [ \SL_2(\Z) : \mathrm{\Gamma_0}(1984) ] = 384$.

\section{Applications of Waldspurger's Theorem}\label{section:appl}

In this section we will present a few examples explaining how to use 
Waldspurger's Theorem.
The idea of using Waldspurger's Theorem for an elliptic curve is motivated by Tunnell's 
famous work on the congruent number problem. 
We will see however that our case needs many 
more computations to get any desired result. In the examples that follow we will first use our algorithm 
in Theorem~\ref{thm:algo} to compute the space of cusp forms that are Shimura equivalent to the given 
elliptic curve and then use Waldspurger's Theorem to get some interesting results. We will follow the notation 
adopted in the previous section. 

\subsection{A First Example}
Our first example will be the elliptic curve $E$ over $\Q$ given by 
\[E: Y^2 = X^3 + X + 1.
 \]
The conductor of $E$ is $496 = 16\times 31$ and $E$ does not have 
complex multiplication.
Let $\phi \in S_2^{\mathrm{new}}(496,\chi_{\mathrm{triv}})$ be the corresponding newform
given by the Modularity Theorem; $\phi$ has the following $q$-expansion:
\[ \phi = q - 3q^5 + 3q^7 - 3q^9 - 2q^{11} - 4q^{13} - q^{19} + O(q^{20}).
 \]
It is to be noted that $\phi$ satisfies the hypothesis $\mathrm{(H1)}$--this 
follows by  
Theorem~\ref{thm:Vig}, and since $16 \mid M_{\phi}$, $\phi$ satisfies $\mathrm{(H2)}$. Let $\chi$ be a Dirichlet 
character with $\chi^2 = \chi_{\mathrm{triv}}$. By Theorem~\ref{thm:Flicker} there exists $N$ such 
that $S_{3/2}(N,\chi,\phi) \ne \{0\}$. Note that we must have $496 \mid (N/2)$.
 
In order to apply Waldspurger's Theorem we would like to compute an eigenbasis for the
summand $S_{3/2}(N,\chi,\phi)$ for a suitable $N$ and $\chi$. We will assume 
$\chi$ to be the trivial character $\chi_{\mathrm{triv}}$. We use Theorem~\ref{thm:algo}
to find out that $S_{3/2}(992,\chi,\phi)= \{0\}$. However at level $1984$ we get that the space 
$S_{3/2}(1984,\chi,\phi)$ has a basis $\{f_1,f_2,f_3\}$ where $f_1$, $f_2$ and $f_3$ have the 
following $q$-expansions:
\begin{equation*}
\begin{split}
f_1 & = q^3 + q^{43} - 2q^{75} + 2q^{83} + q^{91} + 3q^{115} - 3q^{123} + O(q^{145}):= \sum_{n=1}^{\infty} a_nq^n,\\
f_2 & = q^{15} + q^{23} - q^{31} + 2q^{55} + q^{79} - 3q^{119} + O(q^{145}):= \sum_{n=1}^{\infty} b_nq^n,\\
f_3 & = q^{17} + q^{57} + q^{65} + 2q^{73} - q^{89} - q^{105} + q^{137} + O(q^{145}):= \sum_{n=1}^{\infty} c_nq^n.
\end{split}
\end{equation*}

We note that the space $S_{3/2}(1984,\chi)$ is $119$-dimensional.

By Waldspurger's Theorem~\ref{thm:wald} there exists a
function $A_{\phi}$ on square-free positive integers $n$ such that 
\[ {A_{\phi}(n)}^{2}=\mathrm{L}(E_{-n},1)
 \]
and
\[ S_{3/2}(1984,\chi,\phi)=\bigoplus \overline{\mathrm{U}}(E,\phi,A_{\phi}),
 \]
where the sum is over all $E \geq 1$ such that $\widetilde{N_{\phi}} \mid E \mid 1984$.
We already know the left-hand side of the above identity. Henceforth we will be 
interested in computing the right-hand side. We will first compute $\widetilde{N_{\phi}}$ and 
then $\overline{\mathrm{U}}(E,\phi,A_{\phi})$ for $\widetilde{N_{\phi}} \mid E \mid 1984$.

We need to compute local components $\widetilde{n_p}$ for each prime $p$. We consider the following cases. 
\begin{enumerate}[leftmargin=1.5cm]
\item[Case 1.] $p$ odd and $p \ne 31$.

In this case $m_p = 0$ and since $p \nmid N$ the local character $\chi_{0,p}$ is unramified. Hence we get 
that $\widetilde{n_p} = 0$.

\item[Case 2.] $p = 31$.

Here $m_{31} =1$. Since $\lambda_{31} \ne 0$ using Corollary \ref{cor:special} it follows that the local
component $\rho_{31}$ is a special representation of $\GL_2(\Q_{31})$ 
and so $31 \in S$. Also, note that $\Zs_{31}/{\Zs_{31}}^2$ is generated by $11$ mod ${\Zs_{31}}^2$ and using 
Proposition~\ref{prop:ev} we can show that $\chi_{0,31}(11) =1$. Thus $\chi_{0,31}$ is unramified and so, 
$\widetilde{n_{31}} = 1$.

\item[Case 3.] $p = 2$.

In this case $m_{2} =4$ and it is clear from the q-expansion of $\phi$ that $\lambda_2 = 0$. 
We need some information about the set $\Omega_2(\phi)$ (see Equation~\ref{eq:Omega_p}).
In our case, looking at $f_1,f_2$ and $f_3$, we get that $\{ 1,3,7\} \subseteq \Omega_2(\phi)$. 
Since $\nu_2(1)=\nu_2(3)=\nu_2(7) = 0$, we get $\widetilde{n_{2}} = m_2 +2 = 6$.
\end{enumerate}

Hence
\[\widetilde{N_{\phi}} = 31\times 2^6 = 1984.\] 

Thus we have $E = \widetilde{N_{\phi}} = 1984$ and we 
would like to know how the space $\overline{\mathrm{U}}(1984,\phi,A_{\phi})$ looks. For that the next 
immediate task will be to compute $\mathrm{U}_p(e_p,\phi)$ where $e_p = \nu_p(1984)$. We consider the 
following cases:

\begin{enumerate}[leftmargin=1.5cm]
\item[Case 1.] $p$ odd and $p \ne 31$.

Here, $e_p = 0$ and $\mathrm{U}_p(0,\phi)$ consists of only one function $c_{p,\lambda_p^{\prime}}^{(0)}$ defined on 
$\Qs_p$. 
Recall that $c_{p,\lambda_p^{\prime}}^{(0)}(n) = 1$ for $n$ square-free.

\item[Case 2.] $p =31$.

In this case $e_{31} = 1$ and as already seen, $31 \in S$ and $\chi_{0,31}$ is unramified. 
So, $\mathrm{U}_{31}(1,\phi) = \{ c_{31,\lambda_{31}^{\prime}}^{(5)} \}$.
Note that $\lambda_{31} = -1$ and hence $\lambda_{31}^{\prime} = 
(31)^{-1/2} \lambda_{31} =-(31)^{-1/2} $. Again using Proposition~\ref{prop:ev} we can show that
$\chi_{0,31}({31}^{-1})= -1$. Also note that $(31,n)_{31} = \kro{n}{31}$. So for $n$ square-free we have,
\begin{eqnarray*}
c_{31,\lambda_p^{\prime}}^{(5)}(n)=\left\{ \begin{array}{ll}
2^{1/2} & \ \ \text{if}\ \nu_{31}(n)=0\ \text{and}\ \kro{n}{31} = -1 \\
0 & \ \ \text{if}\ \nu_{31}(n)=0\ \text{and}\ \kro{n}{31} = 1\\
1 & \ \ \text{if}\ \nu_{31}(n)=1. \\
\end{array}\right.
\end{eqnarray*}

\item[Case 3.] $p=2$.

Here $e_{2} = 6$. Since $\lambda_2 =0$ and $\{ 1,3,7\} \subseteq \Omega_2(\phi)$, we see that
$\mathrm{U}_{2}(6,\phi)$ consists of $\gamma_{0,1}, \gamma_{0,3}, \gamma_{0,7}$ which are
the characteristic functions of residue classes of $1,3,7$ modulo $8$ respectively.
By our methods so far we do not know whether $5$ belongs to $\Omega_2(\phi)$ or not.
\end{enumerate}

Recall that $\overline{\mathrm{U}}(E,\phi,A_{\phi})$ is the space generated by the functions 
$f(\underline{c},A_{\phi})$ where $\underline{c} \in \prod_p{\mathrm{U}_p(e_p,\phi)}$.
Thus in our case $\underline{c} = (c_p)_p$ where, for odd primes $p \ne 31$ we have
$c_p= c_{p,\lambda_p^{\prime}}^{(0)}$, $c_{31} = c_{31,\lambda_{31}^{\prime}}^{(5)}$ and
for $c_2$ the possible choices are $\gamma_{0,1}$, $\gamma_{0,3}$, $\gamma_{0,5}$ and $\gamma_{0,7}$.
By using Waldspurger's Theorem \ref{thm:wald} we have
\[S_{3/2}(1984,\chi,\phi)=\overline{\mathrm{U}}(1984,\phi,A_{\phi})
 \]
and so every cusp form in the space on the left-hand side  
can be written in terms of
\[
f(\underline{c},A_{\phi})(z) :=
\sum_{n=1}^{\infty}{ A_{\phi}(n^{\mathrm{sc}})n^{1/4}\prod_p{c_p(n)}\ q^n}\]
for some $\underline{c} = (c_p) \in \prod U_p(e_p,\phi)$.

We use Theorem~\ref{thm:congruence} to conclude that $f_1$ has non-zero $n$-th 
coefficients only 
for $n \equiv 3 \pmod{8}$, $f_2$ has non-zero coefficients only for $n \equiv 7 \pmod{8}$ and 
$f_3$ has non-zero coefficients only for $n \equiv 1 \pmod{8}$. 

Since $f_1$ has non-zero $a_n$ only for $n \equiv 3 \pmod{8}$, taking 
$\underline{c}$ as above with $c_2 = \gamma_{0,3}$ we get that
for $n$ square-free,
\begin{equation}\label{eq:an}
\begin{split}
a_n = \beta_1A_{\phi}(n)n^{1/4} & c_2(n)c_{31}(n) =\\ 
&\begin{cases}
2^{1/2}\beta_1A_{\phi}(n)n^{1/4} & \ \text{if}\ \nu_{31}(n)=0,\ \kro{n}{31}=-1\ \text{and}\ n\equiv3\pmod{8}\\
\beta_1A_{\phi}(n)n^{1/4} & \ \text{if}\ \nu_{31}(n)=1\ \text{and}\ n\equiv3\pmod{8}\\
0 & \text{otherwise},
\end{cases}
\end{split}
\end{equation}
for some complex constant $\beta_1$.
Similarly, taking $c_2 = \gamma_{0,7}$ for $f_2$ and $c_2 = \gamma_{0,1}$ for $f_3$ respectively 
we get that
\begin{equation}\label{eq:bn}
\begin{split}
b_n = \beta_2A_{\phi}(n)n^{1/4} & c_2(n)c_{31}(n) =\\ 
&\begin{cases}
2^{1/2}\beta_2A_{\phi}(n)n^{1/4} &\ \text{if}\ \nu_{31}(n)=0,\ \kro{n}{31}=-1\ \text{and}\ n\equiv7\pmod{8}\\
\beta_2A_{\phi}(n)n^{1/4} &\ \text{if}\ \nu_{31}(n)=1\ \text{and}\ n\equiv7\pmod{8}\\
0 & \text{otherwise},
\end{cases}
\end{split}
\end{equation} 
for some complex constant $\beta_2$ and 
\begin{equation}\label{eq:cn}
\begin{split}
c_n = \beta_3A_{\phi}(n)n^{1/4} & c_2(n)c_{31}(n) =\\ 
&\begin{cases}
2^{1/2}\beta_3A_{\phi}(n)n^{1/4} &\ \text{if}\ \nu_{31}(n)=0,\ \kro{n}{31}=-1\ \text{and}\ n\equiv1\pmod{8}\\
\beta_3A_{\phi}(n)n^{1/4} &\ \text{if}\ \nu_{31}(n)=1\ \text{and}\ n\equiv1\pmod{8}\\
0 & \text{otherwise},
\end{cases}
\end{split}
\end{equation}
for some complex constant $\beta_3$.

We have the following proposition which allows us to calculate the critical values 
of the $\mathrm{L}$-functions of $E_{-n}$, the $(-n)$-th quadratic twists of $E$.

\begin{prop}\label{prop:waldapp1}
Let $E$ be as above and $n$ be a positive square-free integer. 
\begin{enumerate}
\item[(i)] If $\nu_{31}(n)=0$, $n\equiv3\pmod{8}$ and $\kro{n}{31} =-1$ then,
       \[\mathrm{L}(E_{-n},1) = \frac{a_n^2 }{2 {\beta_1}^2\sqrt{n}}.\]
\item[(ii)] If $\nu_{31}(n)=1$, $n\equiv3\pmod{8}$ then,
        \[\mathrm{L}(E_{-n},1) = \frac{a_n^2}{{\beta_1}^2\sqrt{n}}.\]
\item[(iii)] If $\nu_{31}(n)=0$, $n\equiv7\pmod{8}$ and $\kro{n}{31} =-1$ then,
        \[\mathrm{L}(E_{-n},1) = \frac{b_n^2}{2{\beta_2}^2\sqrt{n}}.\]
\item[(iv)] If $\nu_{31}(n)=1$, $n\equiv7\pmod{8}$ then,
        \[\mathrm{L}(E_{-n},1) = \frac{b_n^2}{{\beta_2}^2\sqrt{n}}.\]
\item[(v)] If $\nu_{31}(n)=0$, $n\equiv1\pmod{8}$ and $\kro{n}{31} =-1$ then,
        \[\mathrm{L}(E_{-n},1) = \frac{c_n^2}{2{\beta_3}^2\sqrt{n}}.\]
\item[(vi)] If $\nu_{31}(n)=1$, $n\equiv1\pmod{8}$ then,
       \[\mathrm{L}(E_{-n},1) = \frac{c_n^2}{{\beta_3}^2\sqrt{n}}.\]
\end{enumerate}
\end{prop}
\begin{proof}
Using Waldspurger's Theorem \ref{thm:wald} we know the existence of a function $A_\phi$ on square-free numbers such that 
${A_{\phi}(n)}^{2}=\mathrm{L}(E_{-n},1)$. The proof follows now using Equations \eqref{eq:an}, \eqref{eq:bn} and 
\eqref{eq:cn}.
\end{proof}

We will show now how we use the above to calculate the order of the Tate-Shafarevich group $\SHA{E_{-n}/\Q}$.
We will be assuming the Birch and Swinnerton-Dyer Conjecture for rank zero elliptic curves:
\begin{equation}\label{eq:BSD}
\mathrm{L}(E_{-n},1) = \frac{\lvert \SHA{E_{-n}/\Q}\rvert \cdot \Omega_{E_{-n}}\cdot \prod_p{c_p}}
{\lvert E_{-n,\mathrm{tor}}\rvert^2}
\end{equation}
where $\Omega_{E_{-n}}$ stands for the real period of $E_{-n}$ (since $E_{-n}(\R)$ is connected), 
$c_p$ for the $p$-th Tamagawa number of $E_{-n}$ and $E_{-n,\mathrm{tor}}$ stands for the torsion group of 
$E_{-n}$, all of which are easily computable.

We have the following lemma.
\begin{lem}\label{lem:tor}
 Let $E:  Y^2 = X^3 +X +1$. Then $E_{n,\mathrm{tor}}=0$ for all 
square-free integers $n$.
\end{lem}
\begin{proof} 
Let $K=\Q(\sqrt{n})$. It is well-known that the map
\[
E_n(\Q) \rightarrow E(K)
\]
given by
\[
O \mapsto O,\qquad (X,Y) \mapsto \left( \frac{X}{n}, \frac{Y}{n \sqrt{n}} \right)
\]
is an injective group homomorphism~\footnote{As the map simply scales the variables,
it takes lines to lines and so must define a homomorphism of Mordell-Weil groups.}. Thus it is sufficient to show that $E(K)$ has
trivial torsion subgroup. Recall that the discriminant of $E$ is $-496=-16 \times 31$.
Let $p \ne 2$, $31$ be a rational prime and let $\mathfrak{P}$ be a prime ideal
of $K$ dividing $p$. Then $E$ has good reduction at $\mathfrak{P}$. Moreover, if 
$e_\mathfrak{P}<p-1$ then the reduction map $E(K)_\mathrm{tor} \rightarrow E(\F_\mathfrak{P})$
is injective \cite[p.501]{Katz}, where $e_\mathfrak{P}$ is the ramification index for $\mathfrak{P}$
and $\F_\mathfrak{P}$ denotes the residue field of $\mathfrak{P}$.
Thus if $p\geq 5$ and $p \ne 31$ then this map is injective. Now we
take $p=5$, $7$, so $E(\F_\mathfrak{P})$ is a subgroup of $E(\F_{25})$
and $E(\F_{49})$ respectively. Using {\tt MAGMA} we find
\[
E(\F_{25}) \cong \Z/3\Z \times \Z/9\Z, \qquad
E(\F_{49}) \cong \Z/55\Z.
\]
Since these two groups have coprime orders, it follows that
$E(K)_\mathrm{tor}=0$ and so $E_{n,\mathrm{tor}}=0$.
\end{proof}

Further, since the discriminant of $E_{-1}$ is $-496=2^4 \times 31$, by Lemma~\ref{lem:celery}
we know that $\Omega(E_{-n})=\Omega(E_{-1})/\sqrt{n}$.

It is clear that the quantity ${\mathrm{L}(E_{-n},1)}/{\Omega_{E_{-n}}}$ 
is an integer. Using {\tt MAGMA} we compute this integer for $n \in \{3,15,17\}$. In particular for $n=3$ we 
get that ${\mathrm{L}(E_{-3},1)}/{\Omega_{E_{-3}}} = 2$. Substituting this in part $(i)$ of Proposition~\ref{prop:waldapp1} 
and using Lemma~\ref{lem:celery}, it follows that $\Omega_{E_{-1}} = 1/4{\beta_1}^2$. Doing similar calculations with 
$n = 15$, $17$ we get
\begin{equation}\label{eq:betas}
 \Omega_{E_{-1}} = \frac{1}{4{\beta_1}^2} = \frac{1}{4{\beta_2}^2} = \frac{1}{8{\beta_3}^2}.
\end{equation}

Now recall that $W(E_{-n}/\Q)$ denotes the root number for elliptic curve $E_{-n}$ over 
rational numbers. We have the following proposition. 
The methods used here to compute the root numbers are well-known and 
we refer to \cite{Connell}.
\begin{prop}\label{prop:root}
For $E$ as above and $n$ positive square-free the following holds.
\begin{enumerate}
\item[(i)] If $\nu_{31}(n)=0$ then,
\[W(E_{-n}/\Q) = \begin{cases}
-1 & n \equiv 1,3,7 \pmod{8},\ \kro{n}{31}=1\ \mathrm{or}\\ 
& n \equiv 5 \pmod{8},\ \kro{n}{31}=-1\ \mathrm{or}\\
& n\ \mathrm{even},\ \kro{n}{31}=-1;\\
1  & n \equiv 1,3,7 \pmod{8},\ \kro{n}{31}=-1\ \mathrm{or}\\ 
& n \equiv 5 \pmod{8},\ \kro{n}{31}=1\ \mathrm{or}\\
& n\ \mathrm{even},\ \kro{n}{31}=1.\\        
\end{cases}
\]
\item[(ii)] If $\nu_{31}(n)=1$ then,
\[W(E_{-n}/\Q) = \begin{cases}
-1 & n \equiv 5 \pmod{8}\ \mathrm{or}\\
& n\ \mathrm{even};\\
1 & n \equiv 1,3,7 \pmod{8}.\\
\end{cases}
\]
\end{enumerate}
\end{prop}

Before computing the order of the Tate-Shafarevich group $\SHA{E_{-n}/\Q}$, we have the following refinement of 
Theorem \ref{prop:waldapp1}.

\begin{thm}\label{thm:waldapp2}
Let $E : Y^2 = X^3 + X +1$
and $f=f_1+f_2+\sqrt{2} f_3=\sum d_n q^n$. 
Then, for positive square-free $n \equiv 1$, $3$, $7 \pmod{8}$
\[
\mathrm{L}(E_{-n},1) = \frac{2^{(\nu_{31}(n)+1)} \Omega_{E_{-1}}}{\sqrt{n}}
\cdot d_n^2.\]
\end{thm}
\begin{proof}
Note that $d_n = a_n +b_n + \sqrt{2}c_n$. 
It is important for the proof to note that $a_n =0$ for $n \not \equiv 3 \pmod{8}$, and $b_n =0$ for $n \not \equiv 7 \pmod{8}$,
and $c_n=0$ for $n \not \equiv 1 \pmod{8}$; we proved this
by applying Theorem~\ref{thm:congruence}.
It follows from equations \eqref{eq:an}, \eqref{eq:bn} and \eqref{eq:cn} that $d_n = 0$ 
whenever $n \equiv 1,3,7 \pmod{8}$ and the Kronecker symbol $\kro{n}{31} = 1$. Further by Proposition \ref{prop:root} if
$n \equiv 1,3,7 \pmod{8}$ and $\kro{n}{31} = 1$ then $W(E_{-n},\Q) = -1$ and so $\mathrm{L}(E_{-n},1) = 0$. 
Thus the theorem follows when $\kro{n}{31} = 1$. 

In the case when $\kro{n}{31} = -1$, the refinement follows by using  
Equation \eqref{eq:betas} in  
Theorem \ref{prop:waldapp1}.
\end{proof}

We have now the following corollary which computes the order of the Tate-Shafarevich group $\SHA{E_{-n}/\Q}$.
\begin{cor}\label{cor:sha}
Let $E : Y^2 = X^3 + X +1$ and $f=f_1+f_2+\sqrt{2} f_3=\sum d_n q^n$. 
Let $n$ be positive square-free number such that $n \equiv 1$, $3$, $7 \pmod{8}$ and $E_{-n}$ has rank zero. Then, 
assuming the Birch and Swinnerton-Dyer conjecture, 
\[
|\SHA{E_{-n}/\Q}|= \frac{2^{(\nu_{31}(n)+1)}}{\prod_{p}{c_p}} \cdot d_n^2
\]
where the Tamagawa numbers $c_p$ of $E_{-n}$ are given by
\[
c_2=\begin{cases} 
1 & n \equiv 3,7 \pmod{8}\\
2 & n \equiv 1,5 \pmod{8},
\end{cases}
\qquad 
c_{31}=\begin{cases}
1 & 31 \nmid n, \\
4 & 31 \mid n, \kro{n/31}{31}=1\\
2 & 31 \mid n, \kro{n/31}{31}=-1,
\end{cases}
\]
and $c_p=\# E_{-1}(\F_p)[2]$ for $p \mid n$, $p \ne 31$, and $c_p=1$ for all other primes $p$.
\end{cor}
\begin{proof}
From Lemma \ref{lem:tor} we have $E_{-n,\mathrm{tor}}=0$ for all square-free integers $n$. Substituting this and 
$\Omega(E_{-n})=\Omega(E_{-1})/\sqrt{n}$ in Equation \eqref{eq:BSD} we get that
\[\lvert \SHA{E_{-n}/\Q}\rvert \ = \ \frac{\mathrm{L}(E_{-n},1) \cdot \sqrt{n}}
{\Omega_{E_{-1}}\cdot \prod_p{c_p}}\  = \  \frac{2^{(\nu_{31}(n)+1)}}{\prod_p{c_p}} \cdot d_n^2 \ ;\]
the last equality follows by Theorem \ref{thm:waldapp2}.

We use Tate's algorithm (see \cite[p.364--p.368]{Silverman2}) 
to compute the Tamagawa numbers $c_p$.
\end{proof}

We have the following easy corollary to Theorem \ref{thm:waldapp2}.

\begin{cor}\label{cor:rank2}
 Suppose $n \equiv 1$, $3$, $7 \pmod{8}$ and $\kro{n}{31} = -1$. Then assuming the Birch and Swinnerton-Dyer Conjecture,
\[\mathrm{Rank}(E_{-n}) \geq 2 \Leftrightarrow d_n = 0.\]
\end{cor}
\begin{proof}
By Proposition \ref{prop:root}, if $n \equiv 1$, $3$, $7 \pmod{8}$ and $\kro{n}{31} = -1$ 
then $W(E_{-n}/\Q) = 1$. Thus the analytic rank is even, and so by BSD,
the rank is even.
 The corollary now follows using 
Theorem \ref{thm:waldapp2}.
\end{proof}

\noindent {\bf Remark.} By Proposition \ref{prop:root} if $n$ is a square-free integer
such that $n \equiv 5 \pmod{8}$ then $\mathrm{L}(E_{-n},1) = 0$ whenever either 
$\nu_{31}(n) =1 $, or $\nu_{31}(n) =0 $ and $\kro{n}{31} = -1$. One can 
also obtain this by using Waldspurger's Theorem. In fact since $f_1,f_2,f_3$ span 
$S_{3/2}(1984,\chi,\phi)$ and none of them have a non-zero coefficient 
for $n \equiv 5 \pmod{8}$ we obtain
\[
 A_{\phi}(n)c_{31}(n) = 0\ \text{whenever}\ n \equiv 5 \pmod{8}.
\]
The statement now follows since $c_{31}(n) \ne 0$ if either $\nu_{31}(n) =1 $, or $\nu_{31}(n) =0 $
and $\kro{n}{31} = -1$. However these methods fail to provide any information~\footnote{In fact doing computations 
using {\tt MAGMA} we get,   
for example that $\mathrm{L}(E_{-n},1) \ne 0$ for $n = 5$, $69$, $101$, $109$, $133$, $157$, $165$; these $n$ 
satisfy the conditions $n\equiv 5 \pmod{8}$ and $\kro{n}{31} = 1$. However for $n = 149$, $173$,
which also satisfy the same two conditions, we get that 
$\mathrm{L}(E_{-n},1)= 0$ (note thus using root number argument $\mathrm{Rank}(E_{-n}) \geq 2$ 
for $n = 149$, $173$). We do not detect a general pattern.} about $\mathrm{L}(E_{-n},1)$ if 
$n \equiv 5 \pmod{8}$ and $\kro{n}{31} = 1$. We hope to predict what happens in these cases by either going 
to higher levels by suitably twisting $E$ or by allowing non-trivial characters. 

We note here that for newforms $\phi$ of weight $k-1$ and odd and square-free level Baruch and Mao~\cite[Theorem 10.1]{B-M} obtain 
Waldspurger-type results for $\mathrm{L}(\phi\otimes \chi_D, \frac{k-1}{2})$ for all fundamental discriminants $D$. In a subsequent 
paper Mao~\cite[Theorem 1.3]{Mao} removes the square-free condition using the generalized Shimura correspondence. 

\subsection{Second Example}

Our second example will be the rational elliptic curve $E$ of conductor $144$ given by 
\[E: Y^2 = X^3 - 1.
 \]
The corresponding newform $\phi$ is given by
\[
\phi = q + 4q^7 + 2q^{13} - 8q^{19} - 5q^{25} + 4q^{31} - 10q^{37} - 8q^{43} + 9q^{49} + O(q^{50}).
\]
Here $M_{\phi}=144$. 
Using Theorem~\ref{thm:algo} for computing Shimura's decomposition, we find that at the level $576$, 
the space $S_{3/2}(576,\chi_{\mathrm{triv}},\phi) \ne \{0\}$; and 
this space has a basis $\{f_1,f_2,f_3,f_4\}$ where $f_1$, $f_2$, $f_3$ and $f_4$
have the following $q$-expansions:
\begin{equation*}
\begin{split}
f_1 & = q-q^{25} + 5q^{49} -6q^{73} -6q^{97} + O(q^{100}):= \sum_{n=1}^{\infty} a_nq^n,\\
f_2 & = q^5 + q^{29} -q^{53} -2q^{77} + O(q^{100}):= \sum_{n=1}^{\infty} b_nq^n,\\
f_3 & = q^{13} -2q^{61} +q^{85} + O(q^{100}):= \sum_{n=1}^{\infty} c_nq^n,\\
f_4 & = q^{17} -q^{41} -q^{89} + O(q^{100}):= \sum_{n=1}^{\infty} d_nq^n.
\end{split}
\end{equation*}

Doing similar calculations as in the previous example we have the following result.

\begin{thm}\label{thm:waldapp3}
Let $E:  Y^2=X^3-1$. Let 
\[f= f_1/\sqrt{6} + f_2 + \sqrt{2}f_3 + \sqrt{3}f_4 := \sum_{n=1}^{\infty}e_nq^n.\] 
Let $n \ne 1$~\footnote{In the case $n =1$ we still have 
$\mathrm{L}(E_{-n},1) = \frac{\Omega_{E_{-1}}}{\sqrt{n}} \cdot e_n^2$, but since
$|E_{-1,\mathrm{tor}}| = 6$ we get that $|\SHA{E_{-n}/\Q}|= \frac{36}{\prod_{p}{c_p}} \cdot e_n^2 $.}
be a positive square-free integer such that $n \equiv 1\ \text{or}\ 2 \pmod{3}$. Then,
\begin{equation}\label{eqn:Len}
\mathrm{L}(E_{-n},1) = \frac{\Omega_{E_{-1}}}{\sqrt{n}} \cdot e_n^2.
\end{equation}

Further assuming BSD, if $E_{-n}$ has rank zero then,
\[
|\SHA{E_{-n}/\Q}|= \frac{4}{\prod_{p}{c_p}} \cdot e_n^2
\]
where the Tamagawa numbers 
$c_2 = 3$ if $ n \equiv 1 \pmod{8}$, $c_2 =1$ if $n \equiv 3,5,7 \pmod{8}$; $c_3 =2$; 
$c_p=\# E_{-1}(\F_p)[2]$ for $p \mid n$, $p \ne 3$; and $c_p=1$ for all other primes $p$.

\end{thm}

\noindent {\bf Remark.} To consider the case when $3 \mid N$, we try to instead work with elliptic curve 
$E_3$. The curve $E_3$ has conductor $36$ and is isogenous to $E_{-1}$. Hence  
$\mathrm{L}(E_{3n},1) = \mathrm{L}(E_{-n},1)$ and $\mathrm{L}(E_{n},1) = \mathrm{L}(E_{-3n},1)$ 
for all positive square-free $n$ coprime to $3$. Thus computation of $\mathrm{L}(E_{-3n},1)$ for all such $n$ 
will lead to a formula for $\mathrm{L}(E_{n},1)$ for all $n$ square-free. Since the hypothesis $\mathrm{(H2)}$ 
is not satisfied we cannot apply Theorem~\ref{thm:wald} to $E_3$. Let $F:=E_3$ and $\phi^\prime$ be the corresponding newform. 
Using Theorem~\ref{thm:algo} we find that $S_{3/2}(72,\chi_{\mathrm{triv}}, \phi^\prime)$ is two-dimensional spanned by $g_1$ and $g_2$ where
\begin{equation*}
\begin{split}
g_1 & = q - 2q^{10} - 2q^{13} + 4q^{22} - q^{25} + 2q^{34} + 4q^{37} + O(q^{40}),\\
g_2 & = q^2 - q^5 - 2q^{14} + q^{17} + 3q^{29} + O(q^{40}). 
\end{split}
\end{equation*}
Let $g = g_1 + g_2 = \sum_{n=1}^{\infty} a_nq^n$. We try to instead apply Corollary~\ref{cor:wald}. Let $I = \{1,5,13,17\}$, then
for each $i$ in $I$ we obtain 
\[\mathrm{L}(F_{-n},1)= \frac{a_n^2 \cdot \mathrm{L}(F_{-i},1)}{a_i^2} \sqrt{\frac{i}{n}} \  \text{ for } n \equiv i \pmod{24}.
\]
Also by root number calculations $\mathrm{L}(F_{-n},1)=0$ for $n \equiv 7,11,13,17 \pmod{24}$. So the interesting cases we are left with 
are $n \equiv j\pmod{24}$ for $j \in J= \{2,10,14,22\}$. We make the following interesting observation. 
Using {\tt MAGMA} for positive square-free $n \le 1000$ we check that up to $30$ decimal places
\[\mathrm{L}(F_{-n},1)= \frac{a_n^2 \cdot \mathrm{L}(F_{-j},1)}{a_j^2} \sqrt{\frac{j}{n}} \  \text{ for } n \equiv j \pmod{24}.
\]
This observation does not follow from Corollary~\ref{cor:wald}, for example 
$\mathrm{L}(F_{-74},1)= 4 \cdot \frac{\mathrm{L}(F_{-2},1)}{\sqrt{37}}$ but $74/2 \notin {\Qs_2}^2$.

\subsection{Example with a Non-Rational Newform}
In this example we start with a non-rational newform $\psi$ and we show that we can get similar formulae as before 
for the critical values of $\mathrm{L}$-functions of $\psi \otimes \chi_{-n}$.

Let $\psi \in S_2^{\mathrm{new}}(62,\chi_{\mathrm{triv}})$ be a newform of weight $2$, level $62$ and trivial character given 
by the following $q$-expansion,
\[
\psi = q - q^2 + aq^3 + q^4 + (-2a + 2)q^5 - aq^6 + 2q^7 - q^8 + (2a - 1)q^9 + O(q^{10})
\]
where $a$ has minimal polynomial $x^2 - 2x - 2$.

As before using Theorem~\ref{thm:algo} we get that 
the space $S_{3/2}(124,\chi_{\mathrm{triv}},\psi) =\langle f \rangle$ where $f$ 
has the following $q$-expansion,
\[
f = q + (a + 1)q^2 - q^4 - 2aq^5 - aq^7 + (-a - 1)q^8 + (a + 1)q^9 - 2q^{10} + 
O(q^{12}).
\]

Note that Waldspurger's theorem is applicable for the newform $\psi$ since $\rho_2$, the local automorphic representation 
of $\psi$ at $2$, is not supercuspidal; this follows since $\nu_2(62)=1$(see Corollary~\ref{cor:supercuspidal}). 

We have the following proposition.
\begin{prop}
 Let $\psi$ and $f:=\sum_{n=1}^{\infty} a_nq^n$ be as above. Let $n$ be square-free such that 
$n \not\equiv 3 \pmod{8}$ and $\kro{n}{31} \ne -1$. Then
\[
\mathrm{L}(\psi \otimes \chi_{-n},1) = \begin{cases}
                                        \frac{\beta}{\sqrt{n}}\cdot a_n^2  & \text{if $\nu_{31}(n)=1$}\\
                                        \frac{\beta}{2 \sqrt{n}}\cdot a_n^2 & \text{if $\nu_{31}(n)=0$}
                                       \end{cases}
\]
where $\beta = 2 \cdot \mathrm{L}(\psi \otimes \chi_{-1},1)$.
\end{prop}

\begin{proof}
The proof follows by calculations similar to those in the previous examples.
\end{proof}

\subsection{Ternary Quadratic Forms and Tunnell-like Formulae}

For a positive-definite integral quadratic form $Q(x_1,\dots,x_m)$ we define its theta series by
\[
\theta_Q(z)= \sum_{n=0}^\infty \# \{ \mathbf{a} \in \Z^m : Q(\mathbf{a})=n \} \cdot q^n; \qquad
q=\exp(2\pi i z).
\]
 Siegel \cite{Siegel} showed that if 
$Q_1$ and $Q_2$ are positive-definite integral ternary quadratic forms both having level $N$, character
$\chi_d$ and belonging to the same genus, then $\theta_{Q_1}-\theta_{Q_2} \in S_{3/2}(N,\chi_d)$. 
Denote by $S_q(N,\chi_d)$
the subspace of $S_{3/2}(N,\chi_d)$ generated by all such differences of theta series.

It is interesting, when applying Waldspurger's Theorem to a weight $2$ cuspform $\phi$ 
to ask whether the relevant modular form
of weight $3/2$ belongs to $S_q(N,\chi_d)$; in this case we would obtain a Tunnell-like
formula expressing the critical values of the $\mathrm{L}$-functions of twists of $\phi$
in terms of ternary quadratic forms. 
We will illustrate this below by presenting several examples. We point out however that 
this is not always possible. 
In particular for the elliptic curve in our first example, 
$E : Y^2 = X^3 + X +1 $, the space $S_{3/2}(1984,\chi_\mathrm{triv},\phi_E)$ has trivial intersection 
with the subspace $S_q(1984, \chi_\mathrm{triv})$. Note that $\mathrm{L}(E,1) = 0$. 
As we mentioned in the Introduction, for elliptic curves of odd and square-free conductor, 
B\"{o}cherer and Schulze-Pillot \cite{B-S2} showed that an inverse Shimura lift comes from ternary quadratic 
forms if and only if the curve has analytic rank zero. In the examples below we consider levels 
that are neither odd and square-free but the result of B\"{o}cherer and Schulze-Pillot still seems to
hold. 

We do not give details of how to compute $S_q(N,\chi_d)$ or the intersection
$S_q(N,\chi_d) \cap S_{3/2}(N,\chi_d,\phi)$. We merely point out that
it is straightforward to compute a basis for the space $S_q(N,\chi_d)$ with the help
of an algorithm of Dickson~\cite{Dickson,Lehman} for computing quadratic forms
of a given level and character up to equivalence. Computing the intersection
with $S_{3/2}(N,\chi_d,\phi)$ is easy using a suitable adaptation of our Theorem~\ref{thm:algo},
and a result of Bungert \cite[Proposition 4]{Bungert} for computing the Hecke action on 
theta series.

We note here that expressing the forms in $S_{3/2}(N,\chi,\phi)$ in terms of ternary quadratic forms 
has a big advantage in running time for the computation of coefficients of such modular forms for large values 
of $n$ and hence for the computation of critical values of the $\mathrm{L}$-functions for large twists of $\phi$.
In Example~\ref{ex:quad576} below the run time for computing the first $10^5$ coefficients of the theta series is 
just $304.200$ seconds on a modest laptop while the same computation takes over $36$ CPU hours if we do not use the 
representation in terms of ternary quadratic forms. Similarly in Example~\ref{ex:level50} the run time for computing 
the first $10^5$ coefficients of the theta series is $358.820$ seconds.

\noindent {\bf Notation.} We will denote by $[a,b,c,r,s,t]$, the ternary quadratic form given 
by $ax^2 +by^2 +cz^2 +ryz+sxz+txy$.

\begin{ex}\label{ex:level50}

Let $E$ be an elliptic curve of conductor $50$ as in Proposition~\ref{prop:50}.
Let $\phi$ be the newform corresponding to $E$,
\[ \phi = q + q^2 - q^3 + q^4 - q^6 - 2q^7 + q^8 - 2q^9 - 3q^{11} + O(q^{12}).
 \]
Note that $\nu_2(50) = 1$  hence $\rho_2$ is not supercuspidal and we can apply Waldspurger's Theorem. 

We get that $\widetilde{N_{\phi}} = 100$ and $S_{3/2}(100,\chi_{\mathrm{triv}},\phi)$ has a basis 
consisting of $f_1$ and $f_2$ where
\[f_1=q + q^4 - q^6 - q^{11} - 2q^{14} + O(q^{15}):= \sum_{n=1}^{\infty} a_nq^n \]
\[f_2=q^2 - q^3 + q^8 - q^{12} + 2q^{13} + O(q^{15}) := \sum_{n=1}^{\infty} b_nq^n. \]

In fact it turns out that $f_1 = (\theta_{Q_1} - \theta_{Q_2})/2$ and $f_2 = (\theta_{Q_3} - \theta_{Q_4})/2$ 
where $Q_i$'s are the ternary quadratic forms in Proposition~\ref{prop:50}. This can now be proved along 
similar lines to Theorem \ref{thm:waldapp2}.

Again we can compute the order of $\SHA{E_{-n}/\Q}$ assuming BSD. For example, we get that
\[|\SHA{E_{-9318}/\Q}| = 33 ^2 = 1089. \]

We can further consider the real quadratic twists $E_n$. For this we work with the elliptic curve 
$E_{-1}$ of conductor 400,
\[E_{-1} : Y^2 = X^3 + X^2 - 48X - 172.\]

We can show that if $5 \nmid n$ then, 
\[
\mathrm{L}(E_{n},1) = \begin{cases}
\frac{\mathrm{L}(E_{1},1)}{\sqrt{n}}\cdot c_n^2   & \kro{n}{5} = 1\\
\mathrm{L}(E_{17},1)\cdot \sqrt{\frac{17}{n}}\cdot c_n^2   & \kro{n}{5} = -1,                
\end{cases}
\]
where $c_n$ is the $n$-th coefficient of the following linear combination of 
theta series of weight $3/2$ and level $1600$ coming from the ternary quadratic forms:
\begin{equation*}
\begin{split}
&
 -\frac{1}{5} \cdot \theta_{[ 5, 5, 17, -2, -4, 0 ]} 
+\frac{1}{5} \cdot \theta_{[ 5, 9, 10, 2, 2, 4 ]} 
+\frac{1}{10} \cdot \theta_{[ 1, 4, 400, 0, 0, 0 ]} 
-\frac{1}{10} \cdot \theta_{[5, 17, 20, -8, 0, -2 ]}\\
&
-\frac{1}{10} \cdot \theta_{[5, 17, 20, 4, 4, 2 ]}
+\frac{1}{10} \cdot \theta_{[ 8, 13, 20, 12, 8, 4 ]}
-\frac{1}{5} \cdot \theta_{[1, 32, 52, -16, 0, 0 ]}
+\frac{1}{5} \cdot \theta_{[ 8, 13, 17, 6, 4, 4 ]} \\
&
+\frac{1}{10} \cdot \theta_{[ 4, 5, 400, 0, 0, -4 ]}
-\frac{1}{10} \cdot \theta_{[ 4, 16, 101, 0, -4, 0 ]}
+\frac{1}{10} \cdot \theta_{[ 400, 100, 1, 0, 0, 0 ]}\\
&
-\frac{1}{10} \cdot \theta_{[ 125, 100, 4, 0, 0, 100 ]}
+\frac{1}{5} \cdot \theta_{[ 89, 56, 9, -4, -2, -44 ]}
-\frac{1}{5} \cdot \theta_{[ 49, 36, 29, 24, 22, 16 ]}\\
&
-\frac{1}{2} \cdot \theta_{[ 400, 13, 8, 4, 0, 0 ]}
-\frac{1}{10} \cdot \theta_{[ 100, 25, 17, 10, 0, 0 ]}
+\frac{1}{10} \cdot \theta_{[ 52, 32, 25, 0, 0, 16 ]}\\
&
+\frac{1}{2} \cdot \theta_{[ 53, 33, 25, -10, -10, -14 ]}
+\frac{1}{2} \cdot \theta_{[ 400, 400, 1, 0, 0, 0 ]}
+\frac{9}{10} \cdot \theta_{[ 400, 25, 16, 0, 0, 0 ]}\\
&
-\frac{1}{2} \cdot \theta_{[ 201, 201, 4, 4, 4, 2 ]}
+\frac{1}{10} \cdot \theta_{[ 224, 89, 9, -2, -8, -88 ]}
-\frac{1}{10} \cdot \theta_{[ 209, 36, 25, 20, 10, 36 ]}\\
&
-\frac{9}{10} \cdot \theta_{[ 129, 100, 16, 0, -16, -100 ]}
-\frac{4}{5} \cdot \theta_{[ 84, 81, 25, 10, 20, 4 ]}
+\frac{4}{5} \cdot \theta_{[ 89, 49, 41, -6, -14, -38 ]}\\ 
&
-\frac{1}{5} \cdot \theta_{[ 400, 29, 16, 16, 0, 0 ]}
+\frac{1}{5} \cdot \theta_{[ 125, 100, 16, 0, 0, 100 ]}
-\frac{2}{5} \cdot \theta_{[ 100, 96, 21, 8, 20, 80 ]}\\
&
+\frac{2}{5} \cdot \theta_{[ 84, 69, 29, 2, 12, 28 ]}
-\frac{2}{5} \cdot \theta_{[ 400, 32, 13, 8, 0, 0 ]}
+\frac{2}{5} \cdot \theta_{[ 117, 52, 32, -16, -24, -44 ]}\\
&
+\frac{1}{5} \cdot \theta_{[ 400, 25, 17, 10, 0, 0 ]}
+\frac{1}{5} \cdot \theta_{[ 212, 48, 17, 8, 4, 48 ]} 
+\frac{1}{10} \cdot \theta_{[ 208, 32, 25, 0, 0, 32 ]}\\
&
-\frac{1}{5} \cdot \theta_{[ 212, 33, 25, -10, -20, -28 ]}
-\frac{1}{10} \cdot \theta_{[ 208, 33, 32, 32, 32, 16 ]}
-\frac{1}{5} \cdot \theta_{[ 113, 52, 32, 16, 8, 52 ]}.
\end{split}
\end{equation*}
Further, using root number arguments, we get that $\mathrm{L}(E_{-5n},1) = 0$ whenever $n \not\equiv 3 \pmod{8}$ and 
$\mathrm{L}(E_{5n},1) = 0$ whenever $n \equiv 5 \pmod{8}$. We have been unable to derive similar formulae for 
$\mathrm{L}(E_{-5n},1)$ when $n \equiv 3 \pmod{8}$ and for $\mathrm{L}(E_{5n},1)$ when $n \not\equiv 5 \pmod{8}$. 
We consider the twist $E_5$ whose conductor is again $50$ and for which the minimum level to obtain non-zero Shimura 
equivalent forms is $500$, however no new information can be obtained from these forms.
\end{ex}

\begin{ex}\label{ex:quad576}
This example formulates Theorem \ref{thm:waldapp3} in terms of ternary quadratic forms. 
Let $E:Y^2=X^3-1$. Let $n$ be a positive square-free integer such that $n \equiv 1\ \text{or}\ 2 \pmod{3}$. Then
\[\mathrm{L}(E_{-n},1) = \frac{\Omega_{E_{-1}}}{\sqrt{n}} \cdot a_n^2\]  
where $a_n$ is the $n$-th coefficient of the cusp form $f$ of weight $3/2$ and level $576$ that 
can be written as follows as a linear combination theta series:
\begin{equation*}
\begin{split}
& f= \sum_{n=1}^{\infty} a_n q^n = \\
& \frac{1}{\sqrt{6}} \cdot \bigg(\frac{1}{2} \cdot \theta_{[ 1, 36, 45, -36, 0, 0 ]}
-\frac{1}{2} \cdot \theta_{[ 4, 9, 37, 0, -4, 0 ]}
+\frac{1}{2} \cdot \theta_{[ 144, 9, 4, 0, 0, 0 ]}
-\frac{1}{2} \cdot \theta_{[ 45, 36, 4, 0, 0, 36 ]}\\
&
-\frac{1}{2} \cdot \theta_{[ 144, 16, 9, 0, 0, 0 ]}
+\frac{1}{2} \cdot \theta_{[ 49, 36, 16, 0, -16, -36 ]}\bigg)
+ \frac{1}{2} \cdot \theta_{[ 144, 29, 5, 2, 0, 0 ]}
-\frac{1}{2} \cdot \theta_{[ 45, 32, 20, -16, -12, -24 ]}\\
&
+ \sqrt{2} \cdot \bigg( \frac{1}{4} \cdot \theta_{[ 144, 13, 13, 10, 0, 0 ]}
-\frac{1}{4} \cdot \theta_{[ 45, 36, 16, 0, 0, 36 ]} \bigg) 
+ \sqrt{3} \cdot \bigg( \frac{1}{6} \cdot \theta_{[ 1, 4, 144, 0, 0, 0 ]}
-\frac{1}{6} \cdot \theta_{[ 4, 4, 37, 0, -4, 0 ]}\\
&
+\frac{1}{6} \cdot \theta_{[ 4, 5, 36, 0, 0, -4 ]}
-\frac{1}{6} \cdot \theta_{[ 4, 13, 13, -10, 0, 0 ]}
+\frac{1}{3} \cdot \theta_{[ 1, 20, 32, -16, 0, 0 ]}
+\frac{1}{6} \cdot \theta_{[ 4, 5, 29, -2, 0, 0 ]}\\
&
-\frac{1}{2} \cdot \theta_{[ 4, 9, 17, -6, 0, 0 ]}
+\frac{1}{6} \cdot \theta_{[ 144, 16, 1, 0, 0, 0 ]}
-\frac{1}{6} \cdot \theta_{[ 16, 16, 9, 0, 0, 0 ]}
-\frac{1}{3} \cdot \theta_{[ 144, 5, 4, 4, 0, 0 ]}\\
&
+\frac{1}{6} \cdot \theta_{[ 37, 16, 4, 0, 4, 0 ]}
+\frac{1}{6} \cdot \theta_{[ 16, 13, 13, 10, 0, 0 ]}
+\frac{1}{6} \cdot \theta_{[ 32, 21, 4, -4, 0, -16 ]}
-\frac{1}{6} \cdot \theta_{[ 29, 16, 5, 0, 2, 0 ]}\\
&
-\frac{1}{2} \cdot \theta_{[ 144, 36, 1, 0, 0, 0 ]}
+\frac{1}{2} \cdot \theta_{[ 144, 9, 4, 0, 0, 0 ]}
-\frac{1}{6} \cdot \theta_{[ 144, 144, 1, 0, 0, 0 ]}
+\frac{1}{6} \cdot \theta_{[ 49, 36, 16, 0, -16, -36 ]}\\
&
+\frac{1}{2} \cdot \theta_{[ 45, 32, 20, -16, -12, -24 ]}
-\frac{1}{2} \cdot \theta_{[ 32, 29, 29, 22, 16, 16 ]}
-\frac{1}{6} \cdot \theta_{[ 80, 32, 9, 0, 0, 32 ]}
+\frac{1}{2} \cdot \theta_{[ 80, 17, 17, -2, -16, -16 ]}\\
&
-\frac{1}{3} \cdot \theta_{[ 41, 32, 20, 16, 20, 8 ]} \bigg).
\end{split}
\end{equation*}
\end{ex}

\begin{ex}
Let $E: Y^2 + Y = X^3 -7$ be an elliptic curve of conductor $27$ and let $\phi$ be the corresponding newform.
Using Corollary \ref{cor:supercuspidal}, we get that the local component of $\phi$ at $2$ is 
not supercuspidal and hence we can apply Waldspurger's Theorem. We have the following proposition.
\begin{prop}
With $E$ as above let $n$ be a square-free integer.
\begin{enumerate}
\item[(i)] Suppose $n \equiv 1 \pmod 3$.
Let $f$ be  given by 
\begin{equation*}
\begin{split}
f= & \sum_{n=1}^{\infty} a_n q^n = 
-\frac{1}{2} \cdot \theta_{[ 1, 6, 15, -6, 0, 0 ]} 
+\frac{1}{2} \cdot \theta_{[ 4, 4, 7, 4, 4, 2 ]} 
+\theta_{[ 27, 27, 1, 0, 0, 0 ]} \\
& -\theta_{[ 28, 27, 4, 0, 4, 0 ]}-\frac{1}{2} \cdot \theta_{[ 27, 7, 4, 2, 0, 0 ]}
-\frac{1}{2} \cdot \theta_{[ 16, 9, 7, -6, -4, -6 ]} 
+\theta_{[ 31, 16, 7, 4, 2, 16 ]}.
\end{split}
\end{equation*}
If either $\nu_2(n)=1$, or $\nu_2(n)=0$ and $n \equiv 1$, $5 \pmod{8}$ then
\[\mathrm{L}(E_{-n},1) = \frac{\mathrm{L}(E_{-1},1)}{\sqrt{n}} \cdot a_n^2.\]
Otherwise, 
\[\mathrm{L}(E_{-n},1) = \frac{\kappa}{\sqrt{n}} \cdot a_n^2\]
where $\kappa =\sqrt{19}\cdot\mathrm{L}(E_{-19},1)$ if $n \equiv 3 \pmod{8}$ and 
$\kappa =\sqrt{7}\cdot\mathrm{L}(E_{-7},1)$ if $n \equiv 7 \pmod{8}$.
\item[(ii)] Suppose $n \equiv 0 \pmod{3}$ and let $n = 3m$. Let $h \in S_{3/2}(324, \chi_{\mathrm{triv}}, \phi)$ be the cusp form having 
the following $q$-expansion 
\[
h = q^3 - q^{21} + 2 q^{30} - q^{39} - 2 q^{48} - q^{57} - 2q^{66} + q^{75} + O(q^{80}) := \sum_{n=1}^{\infty} b_n q^n.
\]
Further suppose $\kro{m}{3} = 1$. 
If either $\nu_2(n)=1$, or $\nu_2(n)=0$ and $n \equiv 1$, $5 \pmod{8}$ then 
\[\mathrm{L}(E_{-n},1) = \mathrm{L}(E_{-21},1) \cdot \sqrt{\frac{21}{n}} \cdot b_n^2.\]
If $n \equiv 3$, $7 \pmod{8}$ then 
\[\mathrm{L}(E_{-n},1) = \frac{\kappa}{\sqrt{n}} \cdot b_n^2\]
where $\kappa =\sqrt{3}\cdot\mathrm{L}(E_{-3},1)$ if $n \equiv 3 \pmod{8}$ and 
$\kappa =\sqrt{39}\cdot\mathrm{L}(E_{-39},1)$ if $n \equiv 7 \pmod{8}$.
\item[(iii)] If $n= 3m$ and $\kro{m}{3} =-1$ then $\mathrm{L}(E_{-n},1) = 0$.
\item[(iv)] If $n \equiv 2 \pmod 3$ then $\mathrm{L}(E_{-n},1) = 0$.
\end{enumerate}
\end{prop}

The proof of (i) and (ii) follows as in the previous examples, while for (iii) and (iv) one can 
use root number arguments. We point out that the cusp form $h$ which appears in (ii) 
does not come from ternary quadratic forms. Moreover since $E$ is isogenous to $E_{-3}$, for 
$n$ positive square-free $\mathrm{L}(E_{n},1) = \mathrm{L}(E_{-3n},1)$. Thus using above proposition 
we are able to compute the critical values $\mathrm{L}(E_{n},1)$ for all $n$ square-free.

\end{ex}

\section{Tables}

In this section we present tables of orders of Tate-Shafarevich group for twists of some of the elliptic curves 
considered in the previous section. We are assuming the BSD conjecture. 

We first consider the elliptic curve $E: Y^2 = X^3+ X+ 1$ and use the formula in Corollary~\ref{cor:sha} 
to obtain the following table of orders of $\Sha(E_{-n}/\Q)$ for positive square-free $n \le 2000$ with 
$n \equiv 1,3,7 \pmod{8}$ and $\mathrm{L}(E_{-n},1) \ne 0$.

\

\noindent $|\SHA{E_{-n}/\Q}| = 1$ for $n =$  15, 17, 23, 31, 43, 57, 65, 79, 89, 91, 105, 137, 145, 151, 155, 161, 179, 201, 
215, 217, 239, 251, 263, 303, 305, 313, 321, 323, 337, 339, 393, 395, 399, 401, 403, 409, 465, 471, 527, 551, 
571, 595, 601, 611, 619, 633, 651, 663, 673, 681, 697, 699, 705, 755, 759, 767, 787, 843, 849, 871, 879, 895, 
921, 953, 959, 991, 1015, 1019, 1057, 1119, 1129, 1153, 1159, 1171, 1193, 1201, 1209, 1235, 1255, 1257, 1271, 
1329, 1339, 1355, 1367, 1385, 1401, 1441, 1479, 1481, 1545, 1553, 1615, 1633, 1641, 1649, 1673, 1689, 1691, 1729, 
1731, 1735, 1751, 1759, 1767, 1779, 1841, 1851, 1887, 1891, 1921, 1939, 1951, 1959

\

\noindent $|\SHA{E_{-n}/\Q}| = 4$ for $n =$  55, 73, 83, 167, 209, 223, 241, 259, 265, 331, 371, 385, 415, 449, 457, 491, 499, 
587, 649, 695, 703, 761, 881, 983, 1023, 1047, 1049, 1067, 1115, 1139, 1145, 1199, 1295, 1297, 1379, 1407, 1439, 
1463, 1483, 1577, 1579, 1603, 1655, 1687, 1703, 1705, 1793, 1811, 1889, 1903, 1913, 1915, 1937, 1979, 1999 

\

\noindent $|\SHA{E_{-n}/\Q}| = 9$ for $n =$  115, 119, 123, 177, 203, 247, 271, 291, 347, 427, 433, 447, 455, 489, 523, 579,
615, 713, 719, 739, 743, 771, 817, 823, 827, 863, 899, 905, 911, 943, 951, 1003, 1065, 1191, 1195, 1231, 1239, 
1267, 1313, 1319, 1391, 1417, 1491, 1505, 1511, 1515, 1531, 1635, 1695, 1711, 1819, 1897, 1977, 1983

\

\noindent $|\SHA{E_{-n}/\Q}| = 16$ for $n =$  353, 463, 643, 647, 859, 947, 1097, 1111, 1147, 1243, 1345,
1363, 1387, 1393, 1419, 1447, 1487, 1571, 1643, 1667, 1697, 1835, 1855, 1943, 1945, 1987

\

\noindent $|\SHA{E_{-n}/\Q}| = 25$ for $n =$  327, 487, 553, 623, 923, 1207, 1263, 1315, 1455, 1543, 1607, 
1627, 1747, 1763, 1995

\

The following is data for higher orders of $\Sha(E_{-n}/\Q)$ which are attained for $n$ as above with 
$n \le 10000$: 

\

\noindent $|\SHA{E_{-n}/\Q}| = 36$ for $n =$  383, 635, 967, 2347, 2351, 2383, 2411, 2563, 3155, 3391, 3743, 
4091, 4367, 4487, 4519, 4591, 4609, 5323, 5327, 5393, 5423, 5467, 5479, 5555, 5657, 5659, 5803, 
5883, 5963, 6691, 6863, 7159, 7215, 7297, 7307, 7343, 7559, 7567, 7607, 7639, 7895, 7963, 8159, 
8283, 8515, 8635, 8683, 9047, 9385, 9631, 9665, 9667, 9787, 9791

\

\noindent $|\SHA{E_{-n}/\Q}| = 49$ for $n =$  1623, 1753, 2337, 2603, 2927, 2999, 3153, 3279, 3347, 3563, 4043, 
4115, 4331, 4507, 4555, 4955, 4971, 5199, 5347, 5595, 5795, 5955, 6131, 6227, 6447, 6593, 6663, 6695, 
7123, 7283, 7545, 7591, 7687, 7951, 8071, 8135, 8259, 8407, 8431, 8455, 8567, 8755, 8835, 8897, 8923, 
9609, 9771, 9827, 9839

\

\noindent $|\SHA{E_{-n}/\Q}| = 64$ for $n =$  1007, 1727, 2183, 2243, 2455, 2555, 2723, 3763, 3905, 4963, 
5137, 5417, 6587, 6935, 7467, 7483, 7811, 8273, 8737, 9343, 9923

\

\noindent $|\SHA{E_{-n}/\Q}| = 81$ for $n =$  1567, 2683, 2931, 3247, 3323, 3547, 3587, 3855, 3867, 6087, 6305, 
6403, 7153, 7223, 7339, 7833, 7993, 8227, 8447, 8779, 8887, 8895, 9327, 9393, 9931

\

\noindent $|\SHA{E_{-n}/\Q}| = 100$ for $n =$  2827, 3463, 4103, 4543, 5207, 5663, 6847, 7415, 8011, 8015, 8335, 8393, 9143, 9323, 9379

\

\noindent $|\SHA{E_{-n}/\Q}| = 121$ for $n =$  2743, 5703, 7451, 7703, 7873, 7903, 8795, 8983, 9755, 9763

\

\noindent $|\SHA{E_{-n}/\Q}| = 144$ for $n =$  3307, 4643, 9497, 9995

\

\noindent $|\SHA{E_{-n}/\Q}| = 169$ for $n =$  3687, 6527

\

\noindent $|\SHA{E_{-n}/\Q}| = 196$ for $n =$  7867, 9355

\

\noindent $|\SHA{E_{-n}/\Q}| = 225$ for $n =$  7143

\

\noindent $|\SHA{E_{-n}/\Q}| = 289$ for $n =$  5003, 6823, 8903

\

Moreover using Corollary~\ref{cor:rank2} we obtain the following list of positive square free $n \le 10000$ with 
$n \equiv 1,3,7 \pmod{8}$ and $\kro{n}{31} = -1$
such that $d_n = 0$. Hence for these values of $n$, $\mathrm{Rank}(E_{-n}) \geq 2$.

\

\noindent 11, 127, 139, 185, 199, 367, 451, 511, 519, 561, 569, 631, 641, 737, 799, 809,
835, 883, 889, 897, 929, 985, 987, 995, 1009, 1081, 1091, 1131, 1137, 1169, 
1177, 1283, 1443, 1499, 1561, 1563, 1639, 1739, 1801, 1871, 1873, 1883, 2207, 
2409, 2441, 2479, 2495, 2571, 2627, 2785, 2905, 2935, 3081, 3121, 3143, 3289, 
3343, 3377, 3431, 3487, 3499, 3551, 3561, 3799, 3927, 3929, 3959, 4145, 4177, 
4209, 4339, 4355, 4395, 4415, 4463, 4481, 4663, 4735, 4811, 4921, 5017, 5169, 
5335, 5345, 5449, 5561, 5579, 5665, 5671, 5779, 5793, 5849, 5889, 5919, 5951, 
5969, 5979, 5995, 6007, 6031, 6153, 6193, 6211, 6289, 6409, 6465, 6491, 6505, 
6719, 6739, 6761, 6857, 6895, 6911, 6959, 6967, 6999, 7023, 7195, 7207, 7265, 
7315, 7331, 7359, 7513, 7601, 7643, 7711, 7777, 7815, 8139, 8201, 8241, 8249, 
8363, 8369, 8507, 8691, 8769, 8807, 8889, 9127, 9129, 9281, 9311, 9313, 9415, 
9417, 9515, 9543, 9551, 9591, 9647, 9795, 9851, 9895

\

Next we consider the curve $E: Y^2= X^3-1$. We use Theorem~\ref{thm:waldapp3} to obtain the 
orders of $\Sha(E_{-n}/\Q)$ for $n \le 100000$ positive square-free with 
$n \equiv 1,2 \pmod{3}$ and $\mathrm{L}(E_{-n},1) \ne 0$. Here we present a table for values of such $n$ 
with $|\SHA{E_{-n}/\Q}| \ge 256$.

\
 
\noindent $|\SHA{E_{-n}/\Q}| = 256$ for $n$ =  33997, 35341, 38821, 48109, 50893, 62261, 62821, 65285, 70573, 71501, 73309, 75493, 77773, 77797, 84157, 85277, 85333, 89045, 90037, 94813, 96613, 97205

\
 
\noindent $|\SHA{E_{-n}/\Q}| = 289$ for $n$ =  12893, 14717, 14845, 27893, 28661, 30029, 37589, 37621, 39821, 41189, 44789, 45293, 45677, 45869, 53149, 53437, 55061, 55313, 58757, 62989, 68141, 68501, 72077, 72301, 72341, 73421, 80317, 80533, 80813, 82141, 85165, 86357, 87485, 87797, 89501, 89909, 93497, 93565, 95021, 95717, 96221, 96989, 97397

\
 
\noindent $|\SHA{E_{-n}/\Q}| = 324$ for $n$ =  34501, 64237, 79693, 82549

\
 
\noindent $|\SHA{E_{-n}/\Q}| = 361$ for $n$ =  18773, 30341, 31541, 31765, 40949, 43517, 43853, 48341, 49789, 58733, 59021, 61949, 63773, 69541, 71693, 75269, 75949, 76957, 78893, 83093, 83597, 86077, 86341, 86813, 86981, 88045, 92357, 93629, 95429, 95957, 96157, 98269

\
 
\noindent $|\SHA{E_{-n}/\Q}| = 400$ for $n$ =  52261, 64693, 66373, 80029

\
 
\noindent $|\SHA{E_{-n}/\Q}| = 441$ for $n$ =  15629, 23957, 24533, 49157, 53549, 66029, 68813, 70853, 71893, 82333, 82781, 86837

\
 
\noindent $|\SHA{E_{-n}/\Q}| = 484$ for $n$ =  83677, 92797

\
 
\noindent $|\SHA{E_{-n}/\Q}| = 529$ for $n$ =  40829, 51869, 70157, 70877, 73517, 76541, 77213, 79901, 83117, 86117

\
 
\noindent $|\SHA{E_{-n}/\Q}| = 576$ for $n$ =  60037, 85669, 99109, 99469

\
 
\noindent $|\SHA{E_{-n}/\Q}| = 625$ for $n$ =  56605, 57221, 60101, 61757, 85853, 92237, 95653

\
 
\noindent $|\SHA{E_{-n}/\Q}| = 729$ for $n$ =  57557, 65309, 69221, 71741, 71837, 82613, 88661, 98573

\
 
\noindent $|\SHA{E_{-n}/\Q}| = 841$ for $n$ =  76733

\
 
\noindent $|\SHA{E_{-n}/\Q}| = 1089$ for $n$ =  74933

\
 
\noindent $|\SHA{E_{-n}/\Q}| = 1225$ for $n$ =  78797


\begin{thebibliography}{}

\bibitem{Atkin-Li} A.\ O.\ L.\ Atkin and W.\ Li, 
{Twists of newforms and pseudo-eigenvalues of W-Operators},
Invent.\ Math.\ {\bf 48} (1978), 221--243.


\bibitem{B-M} E.\ Baruch and Z.\ Mao, 
{Central value of automorphic L-functions}, 
Geom.\ Funct.\ Anal.\ {\bf 17} (2007), 333--384.


\bibitem{Basmaji} J.\ Basmaji,
{\em Ein Algorithmus zur Berechnung von Hecke-Operatoren und Anwendungen 
auf modulare Kurven},
Ph.\ D. Dissertation, Universit\"{a}t Gesamthochschule Essen, M\"{a}rz 1996.


\bibitem{B-S2} Siegfried B\"{o}cherer and Rainer Schulze-Pillot,
{\em Vector valued theta series and Waldspurger’s theorem}, 
Abh.\ Math.\ Semin.\ Univ.\ Hambg.\ {\bf 64} (1994), 211--233.


 
\bibitem{MAGMA} W.\ Bosma, J.\ Cannon and C.\ Playoust, 
{\em The Magma Algebra System I: The User Language}, 
J.\ Symb.\ Comp.\ {\bf 24} (1997), 235--265. 
(See also {\tt http://magma.maths.usyd.edu.au/magma/})

\bibitem{Bump} D.\ Bump,
{\em Automorphic Forms and Representations},
Cambridge Studies in Advanced Mathematics {\bf 55}, Cambridge University
Press, 1996.

\bibitem{Bungert} M.\ Bungert,
{\em Construction of a cuspform of weight 3/2},
Arch.\ Math.\ {\bf 60} (1993), 530--534.


\bibitem{Cohn} H.\ Cohn,
{A classical invitation to algebraic numbers and class fields},
Springer-Verlag, 1980.

\bibitem{Connell} I.\ Connell,
{Calculating root numbers of elliptic curves over $\Q$},
Manuscripta Math.\ {\bf 82} (1994), 93--104.






\bibitem{Dickson} L.\ E.\ Dickson,
{\em Studies in the theory of numbers},
The University of Chicago Press, Chicago, 1930.



\bibitem{Flicker} Y.\ Flicker, 
{\em Automorphic forms on covering groups of GL(2)}, 
Invent.\ Math.\ {\bf 57} (1980), 119--182.


\bibitem{Hamieh} A.\ Hamieh,
{\em Ternary quadratic forms and half-integral weight modular forms}, 
LMS J.\ Comp.\ Math.\ {\bf 15} (2012), 418--435.


\bibitem{Katz} N.\ Katz,
{\em Galois properties of torsion points
on Abelian varieties}, 
Invent.\ Math.\ {\bf 62} (1981), 481--502.

\bibitem{Kob} N.\ Koblitz,
{\em Introduction to Elliptic Curves and Modular Forms},
GTM {\bf 97}, Springer-Verlag, 1993.

\bibitem{Kohnen} W.\ Kohnen, 
{\em Newforms of half-integral weight}, 
J.\ Reine Angew.\ Math.\ {\bf 333} (1982), 32--72.



\bibitem{Lehman} J.\ Larry Lehman,
{\em Levels of positive definite ternary quadratic forms},
Math.\ Comp.\ {\bf 58} (1992), 399-417.

\bibitem{Mao} Z.\ Mao, 
{A generalized Shimura correspondence for newforms}, 
J.\ Number Theory {\bf 128} (2008), 71--95.





\bibitem{SomaI} S.\ Purkait,
{\em On Shimura's decomposition},
To appear in Int.\ J.\ Number Theory.



\bibitem{Siegel} C.\ L.\ Siegel,
{\em \"{U}ber die analytische Theorie der quadratischen Formen},  
Gesammelte Abhandlungen Bd.\ {\bf 1}, Berlin-Heidelberg-New York: Springer, 1966, 326--405. 

\bibitem{Serre} J.\ P.\ Serre,
{\em A Course in Arithmetic}, 
GTM {\bf 7}, Springer-Verlag, 1973.


\bibitem{Shimura} G.\ Shimura,
{\em On Modular forms of half integral weight},
Annals of Math.\ {\bf 97} (1973), 440--481.

\bibitem{ShimuraII} G.\ Shimura,
{\em The critical values of certain zeta functions associated
with modular forms of half-integral weight},
J.\ Math.\ Soc.\ Japan {\bf 33} (1981), pp 649--672.

\bibitem{Silverman} J.\ H.\  Silverman,
{\em The Arithmetic of Elliptic Curves},
GTM {\bf 106}, Springer-Verlag,  1986.

\bibitem{Silverman2} J.\ H.\  Silverman,
{\em Advanced Topics in the Arithmetic of Elliptic Curves},
GTM {\bf 151}, Springer-Verlag,  1994.


\bibitem{Sturm} J.\ Sturm,
{\em On the congruence of modular forms. Number theory (New York, 1984--1985)},
Lecture Notes in Math.\ {\bf 1240}, Springer, Berlin, (1987), 275--280.

\bibitem{CF} J.\ Tate,
{\em Fourier analysis in number fields and Hecke's zeta-functions},
in {\em Algebraic Number Theory} (J.\ W.\ S.\ Cassels, A.\ Fr\"{o}hlich, eds.) 
Academic Press, 1967, 305--347.

\bibitem{Tate} J.\ Tate,
{\em Number theoretic background},
in {\em Automorphic forms, representations, and L-functions},
Proc.\ Symp.\ in Pure Math.\ XXXIII {\bf 2} (1979), 3--26.

\bibitem{Tunnell} J.\ B.\ Tunnell,
{\em A classical Diophantine problem and modular forms
of weight 3/2},
Invent.\ Math.\ {\bf 72} (1983), 323--334.


\bibitem{Vigneras} M.\ F.\ Vigneras,
{\em Valeur au centre de sym\'{e}trie des functions L associ\'{e}es aux
formes modulaires, S\'{e}minaire de Th\'{e}orie de Nombres, Paris, 1979-1980},
Progr.\ Math.\ {\bf 12}, Birkh\"{a}user, Boston, (1981), 331-356.

\bibitem{Waldspurger} J.\ L.\ Waldspurger,
{\em Sur les coefficients de Fourier des formes
modulaires de poids demi-entier},
J.\ Math.\ Pures Appl.\ {\bf 60} (1981), 375--484.

\bibitem{Yo} S.\ Yoshida, 
{\em Some variants of the congruent number probelm II},
Kyushu J.\ Math.\ {\bf 56} (2002), 147--165.

\end{thebibliography}
\end{document}